\documentclass[fontsize=11pt, titlepage=false]{scrartcl}

\usepackage[utf8]{inputenc}
\usepackage[T1]{fontenc}
\usepackage[greek, USenglish]{babel}

\usepackage[blocks]{authblk}

\title{Strictly hyperbolic Cauchy problems with coefficients low-regular in time and space\thanks{\textbf{Mathematics Subject Classification (2012).} 35S05, 35L30, 47G30}~\thanks{\textbf{Keywords.} Cauchy problem, modulus of continuity, Zygmund space, low-regular, strictly hyperbolic, higher order}}
\author{Daniel Lorenz\thanks{\texttt{daniel.lorenz@math.tu-freiberg.de}}}
\affil{\small{TU Bergakademie Freiberg, Faculty of Mathematics and Computer Science\\Institute of Applied Analysis\\09599 Freiberg, Germany}}

\date{}
\setkomafont{title}{\rmfamily}
\setkomafont{section}{\rmfamily\Large}

\usepackage{graphicx}		
\usepackage{mathscinet}		
\usepackage{enumitem}		
\usepackage{color}			
\usepackage{hyperref}		
\hypersetup{
	colorlinks = true,
	allcolors = blue,
}
\usepackage{pdflscape}		

\usepackage{amsfonts,latexsym}
\usepackage{amssymb}
\usepackage[ntheorem]{empheq}
\usepackage[amsmath, thref, hyperref, thmmarks]{ntheorem}

\usepackage{booktabs}
\usepackage{tabulary}
\usepackage{tabularx}
\usepackage{tablefootnote}
\usepackage{longtable}

\usepackage[font=small,floatrowsep=qquad,captionskip=5pt]{floatrow}
\newcommand{\N}{{\mathbb{N}}}  
\newcommand{\Z}{{\mathbb{Z}}}  
\newcommand{\R}{{\mathbb{R}}}  
\newcommand{\Czi}{{C^{\infty}_{0}}}  

\newcommand{\rmd}{\mathrm{d}}
\newcommand{\I}{\mathrm{i}}

\newcommand{\jbl}{\langle}
\newcommand{\jbr}{\rangle}

\newcommand{\ve}{\varepsilon}

\newcommand{\Sy}{\mathcal{S}}

\newcommand{\OPS}{\Psi}
\newcommand{\OP}{Op}
\newcommand{\jxi}{\jbl\xi\jbr}

\newcommand{\Zyg}{C_\ast}
\newcommand{\Zhyp}{Z_{\text{hyp}}}
\newcommand{\Zpd}{Z_{\text{pd}}}

\newcommand{\chixi}{\chi\Big(\frac{t}{N \eta(\jxi^{-1})}\Big)}

\renewcommand{\Re}{\operatorname{Re}}

\DeclareMathOperator*{\supp}{supp}

\newcommand*\xbar[1]{%
	\hbox{%
		\vbox{%
			\hrule height 0.5pt 
			\kern0.5ex
			\hbox{%
				\kern-0.1em
				\ensuremath{#1}%
				\kern-0.1em
			}%
		}%
	}%
}

\theoremstyle{plain}
\theoremheaderfont{\normalfont\bfseries}
\theorembodyfont{\itshape}
\theoremnumbering{arabic}
\theoremsymbol{}
\theoremseparator{.}

\newtheorem{Theorem}{Theorem}[section]
\newtheorem{Corollary}[Theorem]{Corollary}

\newtheorem{Proposition}[Theorem]{Proposition}

\theoremstyle{plain}
\theoremheaderfont{\normalfont\bfseries}
\theorembodyfont{\upshape}
\theoremnumbering{arabic}
\theoremsymbol{}
\theoremseparator{.}

\newtheorem{Definition}[Theorem]{Definition}
\newtheorem{Remark}[Theorem]{Remark}

\theoremstyle{nonumberplain}
\theoremsymbol{\ensuremath{\Box}}
\theoremheaderfont{\scshape}
\newtheorem{Proof}{Proof}

\numberwithin{equation}{section}

\KOMAoptions{
	bibliography=totoc
}

\allowdisplaybreaks

\begin{document}
	\maketitle

	\begin{abstract}
		We consider the strictly hyperbolic Cauchy problem
		\begin{equation*}
		\begin{cases}
		D_t^m u - \sum\limits_{j = 0}^{m-1} \sum\limits_{|\gamma|+j = m} a_{m-j,\,\gamma}(t,\,x) D_x^\gamma D_t^j u = 0,\\
		D_t^{k-1}u(0,\,x) = g_k(x),\,k = 1,\,\ldots,\,m,
		\end{cases}
		\end{equation*}
		for $(t,\,x) \in [0,\,T]\times \R^n$ with coefficients belonging to the Zygmund class $\Zyg^s$ in $x$ and having a modulus of continuity below Lipschitz in $t$. Imposing additional conditions to control oscillations, we obtain a global (on $[0,\,T]$) $L^2$ energy estimate without loss of derivatives for $s \geq \{1+\ve,\,\frac{2m_0}{2-m_0}\}$, where $m_0$ is linked to the modulus of continuity of the coefficients in time.
	\end{abstract}

	\section{Introduction}\label{Intro}
	
	In this paper we study the strictly hyperbolic Cauchy problem
	\begin{equation}\label{Intro:eq1}
	\begin{cases}
	D_t^m u - \sum\limits_{j = 0}^{m-1} \sum\limits_{|\gamma|+j = m} a_{m-j,\,\gamma}(t,\,x) D_x^\gamma D_t^j u = 0,\\
	D_t^{k-1}u(0,\,x) = g_k(x),\,k = 1,\,\ldots,\,m,
	\end{cases}
	\end{equation}
	for $(t,\,x) \in [0,\,T]\times \R^n$. We aim to derive an global (in time) energy estimate for this equation, when the coefficients are not smooth but their regularity with respect to time is below Lipschitz and with respect to $x$ they belong to the Zygmund space $\Zyg^s$. The derived energy estimate implies in a more or less standard way well-posedness in $H^{m-1} \times H^{m-2} \times \ldots\times L^2$.
			
	Historically, many papers have been devoted to study the well-posedness of the second order strictly hyperbolic Cauchy problem
	\begin{equation}\label{CP:intro0}
	u_{tt} - \sum\limits_{k,\,l = 1}^n \partial_{x_l}(a_{kl}(t,\,x) u_{x_k}) = 0,\quad u(0,x) = u_0(x),\,u_t(0,\,x) = u_1(x),\, (t,\,x) \in [0,\,T]\times\R^n.
	\end{equation}
	It is well-known that, if the coefficients $a_{kl}$ are Lipschitz in $t$ and only measurable in $x$, then the Cauchy problem \eqref{CP:intro0} is well-posed in $H^1 \times L^2$ (see e.g. \cite{Hormander.1963, Mizohata.1973, Hurd.1968}). If the $a_{kl}$'s are Lipschitz continuous in time and $B^\infty$ in space (i.e. $C^\infty$ and all derivatives are bounded), then one can prove well-posedness in $H^s\times H^{s-1}$ for all $s \in \R$. Furthermore, in that case one is able to obtain an energy estimate with no loss of derivatives.
	
	If we lower the regularity of the coefficients below Lipschitz in time, then, generally, we have some loss of derivatives.
	
	There have been different approaches to describe the non-Lipschitz behavior of the coefficients. For the $x$-independent case, the authors of \cite{Colombini.1979} proposed the so called LogLip-property, that is, the coefficients $a_{kl}=a_{kl}(t)$ satisfy
	\begin{equation}\label{Intro:cond0}
	|a_{kl}(t+t_0) - a_{kl}(t)| \leq C t_0 |\log(t_0)|,
	\end{equation}
	for all $t_0 >0,\,t,\,t+t_0 \in [0,\,T]$. Concerning coefficients that also depend on $x$, the authors of \cite{Colombini.1995} and \cite{Agliardi.2004b} considered problem \eqref{Intro:eq1} with $B^\infty$ coefficients in $x$ and LogLip-regularity in $t$. They proved well-posedness in Sobolev spaces and $C^\infty$, respectively, without describing the exact loss of derivatives.
	For second order equations like \eqref{Intro:eq1}, the authors of \cite{Cicognani.2006} provided a classification, linking the loss of derivatives to the modulus of continuity of the coefficients with respect to time.
	For coefficients that are low regular in time and space, the authors of \cite{Colombini.1995} considered coefficients isotropic Log-Lipschitz in $(t,\,x)$. They proved an energy estimate with loss of derivatives which forces us to limit $t$ to a small interval when we want to get an existence result.
	More recently, in \cite{Cicognani.2017} the authors studied strictly hyperbolic Cauchy problems of general order $m$ with coefficients low-regular in time (i.e. Lipschitz or less regular) and smooth in space.
	Most importantly, all results that follow this first approach have in common that, in general, there appears a loss of derivatives when the coefficients are less regular than Lipschitz in $t$.
	
	Another approach to describe the low-regularity of the coefficients with respect to time goes back to \cite{Colombini.2002b}. Generalizing their approach to \eqref{Intro:eq1} gives well-posedness if the coefficients satisfy
	\begin{equation*}
	a_{kl} = a_{kl}(t),\quad a_{kl} \in C[0,\,T]\cap C^1(0,\,T],\quad |a^\prime_{kl}(t)| \leq C t^{-1} \text{ for } t\in(0,\,T].
	\end{equation*}
	In general, also that case comes with a loss of derivatives. However, the $C[0,\,T]\cap C^1(0,\,T]$ regularity in $t$ is not sufficient to give a precise knowledge of the loss of derivatives.
	
	If we assume that the coefficients are more regular, for example $a_{kl} \in C[0,\,T]\cap C^2(0,\,T]$, then we are able to link the loss of derivatives to the control of oscillations (\cite{Colombini.2003e, Hirosawa.2003}).
	
	\begin{Definition}\label{Intro:def1}
		Assume that
		\begin{equation}\label{Intro:cond1}
		a_{kl} = a_{kl}(t),\quad a_{kl} \in C[0,\,T]\cap C^2(0,\,T],\quad |a^{(q)}_{kl}(t)| \leq C \Big(t^{-1}\Big(\log\big(t^{-1}\big)\Big)^\gamma\Big)^q,
		\end{equation}
		for $t \in(0,\,T]$ and $q = 1,\,2$. We say that the coefficients $a_{kl}$ oscillate very slowly, slowly, fast or very fast, if $\gamma = 0,\,\gamma \in (0,\,1),\,\gamma = 1$, if \eqref{Intro:cond1} is not satisfied for $\gamma = 1$, respectively.
	\end{Definition}

	\begin{Theorem}[\cite{Reissig.2003}]\label{Intro:Th1}
		Consider problem \eqref{Intro:eq1} with coefficients satisfying \eqref{Intro:cond1}. Then the energy inequality
		\begin{equation*}
		E_{\nu-\nu_0}(u)(t) \leq C(T) E_\nu(u)(0),
		\end{equation*}
		holds for all $t \in (0,\,T]$ and large $\nu$, where
		\begin{equation*}
		E_\nu(u)(t) := \|(\nabla_x u(t,\,\cdot),\,u_t(t,\,\cdot))\|_{H^\nu}.
		\end{equation*}
		Moreover, we have $\nu_0 =0$ if $\gamma = 0$; $\nu_0$ is an arbitrarily small positive constant if $\gamma \in (0,\,1)$; $\nu_0 $ is a positive constant if $\gamma = 1$; there does not exist a positive constant $\nu_0$ if $\gamma > 1$, i.e. the loss of derivatives is infinite.
	\end{Theorem}

	The paper \cite{Hirosawa.2004} answers the question if it is possible to lower the $C^2$ regularity of the coefficients but still keep the characterization of oscillations and its connection to the loss of derivatives. The author proves that instead of working with $C^2$ regularity from \eqref{Intro:cond1} it is sufficient to work with a suitable $C^{1,\,\beta-1}$ regularity, $\beta\in(1,\,2]$, without changing the statements of Definition~\ref{Intro:def1} and Theorem~\ref{Intro:Th1}

	Finally, in \cite{Kinoshita.2005} the authors provide a general $C^{1,\,\eta,\,\mu,\,\varrho}$ theory that combines a global condition like \eqref{Intro:cond0} for a general modulus of continuity, with generalized local conditions like \eqref{Intro:cond1}.

	For our purpose, this second approach to describe non-Lipschitz behavior of the coefficients is more helpful. By adding a control of oscillations, it is possible to obtain global (in time) energy estimates with no loss of derivatives in certain cases, when the coefficients are $B^\infty$ in $x$ (see e.g. \cite{Kinoshita.2005}). 

	The starting point of our investigations in this paper is the work \cite{Kinoshita.2005}; more precisely, we start from Theorem 2.2 in \cite{Kinoshita.2005} which provides conditions under which the strictly hyperbolic Cauchy problem
	\begin{equation}\label{CP:intro1}
	u_{tt} - \sum\limits_{k,\,l = 1}^n a_{kl}(t,\,x) u_{x_k x_l} = 0,\quad u(0,x) = u_0(x),\,u_t(0,\,x) = u_1(x),
	\end{equation}
	is well-posed in $L^2$ (without loss of derivatives) with coefficients smooth in $x$ but less regular than Lipschitz in $t$. To obtain this result the authors of \cite{Kinoshita.2005} imposed additional local conditions on the first derivative with respect to time of the coefficients.
	
	We aim to generalize the results of \cite{Kinoshita.2005} to equations of general order $m$ and with coefficients low-regular in $x$. By this we mean, that we assume that the coefficients belong to the Zygmund space $\Zyg^s$ with respect to $x$ for some $s>0$. Our goal is to make $s$ as small as possible.
	
	\begin{Definition}\label{Def:Zyg}
		The Zygmund space $\Zyg^s,\,s \geq 0$ consists of all functions $u \in C^{[s]}$ such that
		\begin{equation*}
			\|u\|_{\Zyg^s} =  \sum\limits_{|\alpha| \leq [s]} \|D_x^\alpha u\|_{L^\infty} + \sum\limits_{|\alpha| = [s]} \sup\limits_{x\neq y}\frac{|D_x^\alpha u(x) - 2 D^\alpha_x(\frac{x+y}{2}) + D_x^\alpha u(y)|}{|x-y|^{s-[s]}} < +\infty,
		\end{equation*}
		where $[s]$ denotes the largest integer strictly smaller than $s$.
	\end{Definition}

	The result we obtain reveals that for coefficients that are sufficiently regular in time (see Corollaries~\ref{Intro:cor1} and \ref{Intro:cor2}), it is possible to lower $s$ to $1+\ve$ for any $\ve >0$. The minimum regularity in time that we have to suppose for this to be true is Hölder continuous with exponent $\frac{1}{3}$. If the coefficients are less regular in time than that, we have to choose larger values for $s$ for our result to be valid.

	If we apply our results to the second order equation
	\begin{equation}\label{discuss:q1}
	u_{tt} - \sum\limits_{k,\,l = 1}^n a_{kl}(t,\,x) u_{x_k x_l} = 0,\quad u(0,x) = u_0(x),\,u_t(0,\,x) = u_1(x),
	\end{equation}
	with coefficients
	\begin{equation*}
		a_{kl}(t,\,x) \in C^\mu([0,\,T];\,\Zyg^s) \cap C^2((0,\,T];\,\Zyg^s),
	\end{equation*}
	we can prove the following results.
	
	\begin{Corollary}\label{Intro:cor1}
		Consider equation \eqref{discuss:q1} with coefficients $a_{kl}$ having the modulus of continuity
		\begin{equation*}
			\mu(r) = r \log\left(\frac{1}{r}\right),
		\end{equation*}
		and satisfying the oscillation conditions
		\begin{equation*}
			\|\partial_t a_{kl}(t,\,\cdot)\|_{\Zyg^s} \lesssim \frac{e^{ \frac{1}{2t}}}{t},
		\end{equation*}
		and
		\begin{equation*}
			\|\partial_t^2 a_{kl}(t,\,\cdot)\|_{\Zyg^s} \lesssim \frac{e^{ \frac{1}{t}}}{t^2},
		\end{equation*}
		for $t \in (0,\,T]$.
		Then the energy inequality
		\begin{equation*}
		\| (\jbl D_x \jbr u(t,\,\cdot),\,D_t u(t,\,\cdot))^T\|_{L^2}\\\leq C_T \| (\jbl D_x \jbr u_0(\cdot),\,D_t u_1(\cdot))^T\|_{L^2},
		\end{equation*}
		holds for all $t\in (0,\,T]$ provided that the index $s$ of the Zygmund space $\Zyg^s$ satisfies
		\begin{equation*}
		s \geq 1+\ve,
		\end{equation*}
		for any $\ve > 0$.
	\end{Corollary}
	
	\begin{Corollary}\label{Intro:cor2}
		Consider equation \eqref{discuss:q1} with coefficients $a_{kl}$ having the modulus of continuity
		\begin{equation*}
			\mu(r) = r^\alpha,\, \alpha\in(0,\,1),
		\end{equation*}
		and satisfying the oscillation conditions
		\begin{equation*}
			\|\partial_t a_{kl}(t,\,\cdot)\|_{\Zyg^s} \lesssim (1-\alpha)^{-\frac{1}{2}} t^{-\frac{2-\alpha}{2(1-\alpha)}},
		\end{equation*}
		and
		\begin{equation*}
			\|\partial_t^2 a_{kl}(t,\,\cdot)\|_{\Zyg^s} \lesssim (1-\alpha)^{-1} t^{-\frac{2-\alpha}{1-\alpha}},
		\end{equation*}
		for $t \in (0,\,T]$. 
		Then the energy inequality
		\begin{equation*}
		\| (\jbl D_x \jbr u(t,\,\cdot),\,D_t u(t,\,\cdot))^T\|_{L^2}\\\leq C_T \| (\jbl D_x \jbr u_0(\cdot),\,D_t u_1(\cdot))^T\|_{L^2},
		\end{equation*}
		holds for all $t\in (0,\,T]$ provided that the index $s$ of the Zygmund space $\Zyg^s$ satisfies
		\begin{equation*}
		s \geq  \begin{cases} 1+ \ve & \text{ if } \alpha \geq \frac{1}{3};\\
		2 \frac{1-\alpha}{1+\alpha} & \text { if }\alpha < \frac{1}{3}.
		\end{cases}
		\end{equation*}	
	\end{Corollary}

	The paper is organized as follows: Section~\ref{RESULTS} states the main results of this paper. Examples and remarks are discussed in Section~\ref{EXAMPLES}. Section~\ref{DEF} reviews some definitions and provides an introduction to the pseudodifferential calculus used in this paper. Finally, in Section~\ref{PROOF} we proceed to prove the theorem of Section~\ref{RESULTS}. In Section~\ref{CONCLUSION} we present some concluding remarks.
	
	\section{Statement of results}\label{RESULTS}
	
		We consider the strictly hyperbolic Cauchy problem
	\begin{equation}\label{CP:eq1}
	\begin{cases}
	D_t^m u - \sum\limits_{j = 0}^{m-1} \sum\limits_{|\gamma|+j = m} a_{m-j,\,\gamma}(t,\,x) D_x^\gamma D_t^j u = 0,\\
	D_t^{k-1}u(0,\,x) = g_k(x),\,k = 1,\,\ldots,\,m,
	\end{cases}
	\end{equation}
	for $(t,\,x) \in [0,\,T]\times \R^n$.

	We assume that
	\begin{enumerate}[label = (A\arabic*)]
		\item $g_k \in H^{m-k},\,k= 1,\,\ldots,\,m$,
		\item the characteristic roots $\tau_k = \tau_k(t,\,x,\,\xi)$ of 
		\begin{equation*}
			\tau^m - \sum\limits_{j = 0}^{m-1} \sum\limits_{|\gamma|+j = m]} a_{m-j,\,\gamma}(t,\,x) \xi^\gamma \tau^j = 0,
		\end{equation*}
		are real and distinct for all $(t,\,x,\,\xi) \in [0,\,T]\times\R^n_x \times \R^n_\xi\setminus\{0\}$,
		\item there exist auxiliary functions $\eta$ and $\varrho$ defined on $(0,\,r_0]$ ($r_0$ small) satisfying 
		\begin{equation*}
			\eta \in C^\infty(0,\,r_0],\,\eta(0) = 0,\,\eta^\prime > 0,\,\eta^{\prime\prime} < 0,\,|\eta^{(k)}(r)| \leq C_k r^{-(k-1)}\eta^\prime(r),\,k \geq 1,
		\end{equation*}
		and
		\begin{equation*}
			\varrho \in C^\infty(0,\,r_0],\,\varrho(0) = 0,\,\varrho^\prime > 0,\,\varrho^{\prime\prime} \leq 0,\,|\varrho^{(k)}(r)| \leq C_k r^{-(k-1)}\varrho^\prime(r),\,k \geq 1,
		\end{equation*} such that
		\item  \label{ass:f} the function $f(r) := \frac{r}{\eta(r)}$ is increasing on $(0,\,r_0]$ and $\lim\limits_{r\rightarrow 0+} f(r) = 0$,
		\item \label{ass:gl} the coefficients satisfy the global condition
		\begin{equation*}
			\begin{aligned}
			\sup\limits_{t,\,t_0,\,t+t_0 \in [0,\,T]} \frac{\bigg\|\partial_\xi^\alpha \Big(\sum\limits_{j = 0}^{m-1} \sum\limits_{|\gamma|+j = m} (a_{m-j,\,\gamma}(t+t_0,\,\cdot) - a_{m-j,\,\gamma}(t,\,\cdot))\frac{\xi^\gamma}{|\xi|^{m-j}}\Big) \bigg\|_{\Zyg^s}}{ \frac{t_0}{\eta(t_0)}}\\ \leq C_{\alpha,\,s} |\xi|^{-|\alpha|},
			\end{aligned}
		\end{equation*}
		for $\xi \in\R^n\setminus\{0\},\,\alpha \in \N^n$,
		\item \label{ass:l1} the coefficients satisfy the local condition
		\begin{equation*}
			\bigg\|\partial_\xi^\alpha \partial_t \Big(\sum\limits_{j = 0}^{m-1} \sum\limits_{|\gamma|+j = m} a_{m-j,\,\gamma}(t,\,\cdot)\frac{\xi^\gamma}{|\xi|^{m-j}}\Big) \bigg\|_{\Zyg^s}^2 \leq -C^2_{\alpha,\,s} |\xi|^{-2|\alpha|} \frac{\rmd}{\rmd t}\Big(\frac{1}{\eta^{-1}(t)}\Big),
		\end{equation*}
		for $\xi \in\R^n\setminus\{0\},\,\alpha \in \N^n,\, t \in (0,\,T]$,
		\item \label{ass:l2} the coefficients satisfy the additional local condition
		\begin{equation*}
			\begin{aligned}
				\sup\limits_{t_0>0,\,\tau,\,\tau+t_0 \in [t,\,T]}\frac{\bigg\|\partial_\xi^\alpha \Big(\sum\limits_{j = 0}^{m-1} \sum\limits_{|\gamma|+j = m} (\partial_t a_{m-j,\,\gamma}(\tau+t_0,\,\cdot) - \partial_t a_{m-j,\,\gamma}(\tau,\,\cdot))\frac{\xi^\gamma}{|\xi|^{m-j}}\Big) \bigg\|_{\Zyg^s}}{\varrho(t_0) }\\\leq- C_{\alpha,\,s} |\xi|^{-|\alpha|} \frac{\rmd}{\rmd t}\Big(\frac{1}{\varrho(\eta^{-1}(t))}\Big),
			\end{aligned}
		\end{equation*}
		for $\xi \in\R^n\setminus\{0\},\,\alpha \in \N^n$ and $t \in (0,\,T]$.
	\end{enumerate}


	\begin{Theorem}\label{th:th}
		Consider the Cauchy problem \eqref{CP:eq1} under the assumptions (A1)-(A7) and suppose that for some $m_0 \in (0,\,1]$ we have
		\begin{align}
			\label{th:c1}
			\frac{1}{\eta(\jxi^{-1})} &\in \Sy^{m_0}_{1,\,0},\\
			\label{th:c2}
			\jxi^{-1}\frac{\rmd}{\rmd t}\Big(\frac{1}{\eta^{-1}(t-\jxi^{-1})}\Big)&\in \Sy^{m_0}_{1,\,0}, \text{ for } t \geq \eta(\jxi^{-1})\\
			\label{th:c3}
			 \varrho(\jxi^{-1})\frac{\rmd}{\rmd t}\Big(\frac{1}{\varrho(\eta^{-1}(t-\jxi^{-1}))}\Big)&\in\Sy^{m_0}_{1,\,0}, \text{ for } t \geq \eta(\jxi^{-1}).
		\end{align}
		Then the energy inequality
		\begin{equation*}
			\| (\jbl D_x \jbr^{m-1} u(t,\,\cdot),\,\ldots,\,D_t^{m-1}u(t,\,\cdot))^T\|_{L^2}\\\leq C_T \| (\jbl D_x \jbr^{m-1} g_1(\cdot),\,\ldots,\,D_t^{m-1}g_m(\cdot))^T\|_{L^2},
		\end{equation*}
		holds for all $t\in (0,\,T]$ provided that the index $s$ of the Zygmund space $\Zyg^s$ satisfies
		\begin{equation}\label{cond:s}
			s \geq \max\Big\{1 + \ve, \frac{2m_0}{2-m_0}\Big\},
		\end{equation}
		for any $\ve > 0$.
	\end{Theorem}

	\begin{Remark}
		We note that the conditions \eqref{th:c1}-\eqref{th:c3} are always satisfied for $m_0 = 1$ (Lemma 5.7. in \cite{Kinoshita.2005}). This yields the following corollary.
	\end{Remark}
	
	\begin{Corollary}\label{Th:cor}
		The energy estimate in Theorem~\ref{th:th} holds true even if we just assume (A1)-(A7) (and not \eqref{th:c1}-\eqref{th:c3}) provided that the index $s$ of the Zygmund space $\Zyg^s$ satisfies
		\begin{equation*}
		s \geq 2.
		\end{equation*}
	\end{Corollary}

	\begin{Remark}
		We note that the terms in condition \eqref{cond:s} for $s$ are in equilibrium, if
		\begin{equation*}
			\eta(r) = r^{\frac{2}{3}},
		\end{equation*}
		that is, if the coefficients are Hölder continuous in time with exponent $1-\frac{2}{3}  = \frac{1}{3}$.
	\end{Remark}

	\section{Examples and remarks}\label{EXAMPLES}
	
	Let us begin with some remarks and explanations regarding the assumptions of the previous theorem.
	
	\begin{Remark}
		To formulate the assumptions (A5)-(A7) we use the two auxiliary function $\eta$ and $\varrho$. Typical examples of $\eta$ are
		\begin{equation*}
		\eta(r) = (\log^{[\widetilde m]}\tfrac{1}{r})^{-1},\quad \eta(r) = (\log\tfrac{1}{r})^{-\alpha},\,\alpha > 0,\quad \eta(r) = r^\beta,\,\beta\in(0,\,1).
		\end{equation*}
		Typical examples of $\varrho$ are the same as for $\eta$ but $\varrho(r) = r$ is also admissible.
	\end{Remark}
		
	\begin{Remark}
		The global condition \ref{ass:gl} states that the coefficients have the modulus of continuity $\frac{r}{\eta(r)}$ in time. Due to the assumptions on $\eta$ and \ref{ass:f} this modulus of continuity is always weaker than Lipschitz, i.e. the coefficients are less regular than Lipschitz.
	\end{Remark}
	
	\begin{Remark}
		The local condition \ref{ass:l1} is basically a control of oscillations.
	\end{Remark}
	
	\begin{Remark}
		The additional local condition \ref{ass:l2} is there to replace a condition on the second derivatives of the coefficients.
		
		If the coefficients are $C^2$ with respect to time, then we can choose $\varrho(r) = r$ and condition \ref{ass:l2} turns into a condition on the second derivative with respect to time of the coefficients.
		
		If the coefficients are not $C^2$, then we cannot choose $\varrho(r) = r$ since the supremum might not exist. In these cases, condition \ref{ass:l2} states that the coefficients are more regular than just $C^1$ and $\varrho$ describes how much more regular (than $C^1$) the coefficients are.
		
		We note that the right hand side of this condition becomes more restrictive (for small $t$) the closer we come to $C^1$ regularity of the coefficients, i.e. the larger $\varrho(r)$ gets for small $r$.
	\end{Remark}
	
	Let us try to get a feeling for these conditions by looking at some examples. First, we note that the choice of the modulus of continuity, i.e. $\frac{r}{\eta(r)}$, and the choice of $\varrho(r)$ are completely independent of each other. Moreover, their choice does not influence the result of the theorem. Choosing a less regular modulus of continuity automatically yields more restrictive conditions \ref{ass:l1} and \ref{ass:l2}. Similarly, choosing $\varrho$ such that we are closer to $C^1$ regularity gives a more restrictive condition \ref{ass:l2}.
	
	We start by considering the local condition \ref{ass:l1} for different moduli of continuity. The results are shown in Table~\ref{table1}. As expected the right hand side of \ref{ass:l1} gets more restrictive, i.e. the term $\sqrt{-\frac{\rmd}{\rmd t} \frac{1}{\eta^{-1}(t)}}$ is smaller for small $t$, the further away we get from Lipschitz regularity, i.e. $\eta(r) = 1$.
	\begin{table}[h!]
		\begin{tabulary}{\textwidth}{LCR}
			\toprule
			modulus of continuity & $\eta(r)$& $\sqrt{-\frac{\rmd}{\rmd t} \frac{1}{\eta^{-1}(t)}}$\\
			\midrule
			$r$	& $1$& case excluded by assumptions on $\eta$\\[5pt]
			$r \log\left(\frac{1}{r}\right)$	& $\Big(\log\Big(\frac{1}{r}\Big)\Big)^{-1}$&$\displaystyle\frac{e^{ \frac{1}{2t}}}{t}$
			\\[5pt]
			$r^\alpha,\,\alpha\in(0,\,1)$	& $ r^{1-\alpha}$& $\displaystyle (1-\alpha)^{-\frac{1}{2}} t^{-\frac{2-\alpha}{2(1-\alpha)}}$\\[5pt]
			\bottomrule
		\end{tabulary}
		\caption{Some examples of the local condition \ref{ass:l1} for different moduli of continuity.}
		\label{table1}
	\end{table}
	
	Next, we look at condition \ref{ass:l2} for different $\eta$ and $\varrho$. The examples are shown in Table~\ref{table2}, where the resulting therms for the right hand side of \ref{ass:l2} are written in the respective table cells below $\eta(r)$ and to the right of $\varrho(r)$. As expected, we see that condition \ref{ass:l2} gets more restrictive, i.e. the terms in the cells of Table~\ref{table2} are growing slower for $t \rightarrow 0+$, the further away we are from $C^2$ regularity, i.e. $\varrho(r) = r$.
	\begin{table}[htbp]
		\begin{tabulary}{\textwidth}{L|CL}
			\toprule
			$\varrho(r)$~\textbackslash~$\eta(r)$& $\eta(r) = \Big(\log\Big(\frac{1}{r}\Big)\Big)^{-1}$ & $\eta(r) = s^{1-\alpha}$\\
			\midrule
			$\varrho(r) =  \Big(\log\Big(\frac{1}{r}\Big)\Big)^{-1}$&$\displaystyle \frac{1}{t^2}$ &$\displaystyle \frac{1}{1-\alpha} t^{-1}$ \\[5pt]
			$\varrho(r) = r^\beta$ & $\displaystyle \frac{\beta}{t^2} e^{\frac{\beta}{t}}$& $\displaystyle\frac{\beta}{1-\alpha} t^{-\frac{1-\alpha+\beta}{1-\alpha}}$  \\[5pt]
			$\varrho(r) = r$	&$\displaystyle \frac{1}{t^2} e^{\frac{1}{t}}$&$\displaystyle\frac{1}{1-\alpha} t^{-\frac{2-\alpha}{1-\alpha}}$  \\[5pt]
			\bottomrule
		\end{tabulary}
		\caption{Examples of the additional local condition \ref{ass:l2} for different $\eta$ and $\varrho$. The resulting terms for $\displaystyle -\frac{\rmd}{\rmd t} \bigg(\frac{1}{\varrho(\eta^{-1}(t))}\bigg)$ are shown in the respective table cells.}
		\label{table2}
	\end{table}

	Let us now consider the conditions~\eqref{th:c1}-\eqref{th:c3}. These three conditions are used in the proof to apply sharp G{\aa}rding's inequality for possibly smaller values of $s$.  As stated in Corollary~\ref{Th:cor}, it is possible to ignore these conditions but then one may only work with $s \geq 2$. 
	
	Table~\ref{table01} reviews some examples of the weight $\frac{1}{\eta(\jxi^{-1})}$, i.e. condition~\eqref{th:c1}, for different moduli of continuity. 	Similarly, Table~\ref{table02} and Table~\ref{table03} present some examples for  condition~\eqref{th:c2} and  condition~\eqref{th:c3}, respectively, for different moduli of continuity and different $\varrho$.	
	\begin{table}[h!]
		\begin{tabulary}{\textwidth}{LCCR}
			\toprule
			modulus of continuity & $\frac{1}{\eta(\jxi^{-1})}$&symbol class& $m_0$\\
			\midrule
			$r \log\left(\frac{1}{r}\right)$&
			$\log(\jxi)$&
			$\bigcap\limits_{\ve >0}\Sy^{\ve}_{1,\,0}$&
			any $m_0 > 0$
			\\[5pt]
			$r^\alpha,\,\alpha\in(0,\,1)$&
			$\jxi^{1-\alpha}$&
			$\Sy^{1-\alpha}_{1,\,0}$&
			$m_0 = 1-\alpha$
			\\[5pt]
			$\Big(\log\Big(\frac{1}{r}\Big)\Big)^{-\alpha},$&
			$\jxi \Big(\log(\jxi)\Big)^{-\alpha}$&
			$\bigcup\limits_{\ve > 0}\Sy^{1-\ve}_{1,\,0}$&
			$m_0 = 1$
			\\
			$\alpha\in(0,\,\infty)$&
			&
			&\\[5pt]
			\bottomrule
		\end{tabulary}
		\caption{Some examples of the weight $\frac{1}{\eta(\jxi^{-1})}$, i.e. condition~\eqref{th:c1} for different moduli of continuity.}
		\label{table01}
	\end{table}

	\begin{table}[h!]
		\begin{tabulary}{\textwidth}{LCCR}
			\toprule
			modulus of continuity &  $\jxi^{-1}  \Big(-\frac{\rmd}{\rmd t}\Big(\frac{1}{\eta^{-1}(t-\jxi^{-1})}\Big)\Big)$&symbol class&$ m_0$\\
			\midrule
			$r \log\left(\frac{1}{r}\right)$&
			$(\log(\jxi))^2$&
			$\bigcap\limits_{\ve >0}\Sy^{\ve}_{1,\,0}$&
			any $m_0 > 0$
			\\[5pt]
			$r^\alpha,\,\alpha\in(0,\,1)$&
			$\jxi^{1-\alpha}$&
			$\Sy^{1-\alpha}_{1,\,0}$&
			$m_0 = 1-\alpha$
			\\[5pt]
			\bottomrule
		\end{tabulary}
		\caption{Some examples of the weight $\jxi^{-1}  \Big(-\frac{\rmd}{\rmd t}\Big(\frac{1}{\eta^{-1}(t-\jxi^{-1})}\Big)\Big)$, i.e. condition \ref{th:c2}, for different moduli of continuity.}
		\label{table02}
	\end{table}
	\begin{table}[htbp]
		\begin{tabulary}{\textwidth}{L|C L}
			\toprule
			$\varrho(r)$~\textbackslash~$\eta(r)$& $\eta(r) = \Big(\log\Big(\frac{1}{r}\Big)\Big)^{-1}$ & $\eta(r) = r^{1-\alpha}$\\
			\midrule
			$\varrho(r) =  \Big(\log\Big(\frac{1}{r}\Big)\Big)^{-1}$&
			$\log(\jxi)$ &
			$\displaystyle\frac{1}{1-\alpha} \frac{\jxi^{1-\alpha}}{\log(\jxi)}$
			\\[10pt]
			$\varrho(r) = r^\beta$ &
			$\beta (\log(\jxi))^2$&
			$\displaystyle\frac{\beta}{1-\alpha} \jxi^{1-\alpha}$ 
			\\[10pt]
			$\varrho(r) = r$&
			$(\log(\jxi))^2$&
			$\displaystyle\frac{1}{1-\alpha} \jxi^{1-\alpha}$
			\\[5pt]
			\bottomrule
		\end{tabulary}
		\caption{Examples of $\varrho(\jxi^{-1}) \Big(-\frac{\rmd}{\rmd t}\Big(\frac{1}{\varrho(\eta^{-1}(t-\jxi^{-1}))}\Big)\Big)$, i.e. condition \eqref{th:c3}, for different $\eta$ and $\varrho$. The resulting terms for $\varrho(\jxi^{-1}) \Big(-\frac{\rmd}{\rmd t}\Big(\frac{1}{\varrho(\eta^{-1}(t-\jxi^{-1}))}\Big)\Big)$ are shown in the respective table cells.}
		\label{table03}
	\end{table}
	Looking at the examples in Table~\ref{table03} we note that the weight $\varrho(\jxi^{-1}) \Big(-\frac{\rmd}{\rmd t}\Big(\frac{1}{\varrho(\eta^{-1}(t-\jxi^{-1}))}\Big)\Big)$ changes slightly if we change $\varrho$. However, the weight mainly depends on our choice of $\eta$. For fixed $\eta$, we can choose one $m_0$ such that the weight $\varrho(\jxi^{-1}) \Big(-\frac{\rmd}{\rmd t}\Big(\frac{1}{\varrho(\eta^{-1}(t-\jxi^{-1}))}\Big)\Big)$ belongs to $\Sy^{m_0}_{1,\,0}$ for any of the given examples for $\varrho$. We conclude from this observation, that the possible values for $s$ are mainly determined by our choice of the modulus of continuity.
	
	Finally, Table~\ref{table04} summarizes the discussed examples for the different weights and conditions. For a given modulus of continuity Tabel~\ref{table04} shows possible choices of $m_0$ and the resulting possible values for $s$.
	
	\begin{table}[h!]
		\begin{tabulary}{\textwidth}{LCR}
			\toprule
			modulus of continuity & possible $m_0$& possible $s$\\
			\midrule
			$r \log\left(\frac{1}{r}\right)$&
			any $m_0 > 0$&
			$s \geq 1+\ve$
			\\[15pt]
			$r^\alpha,\,\alpha\in(0,\,1)$&
			$m_0 = 1 - \alpha$&
			$s \geq \max\{1+\ve,\,2\frac{1-\alpha}{1+\alpha}\}$
			\\[15pt]
			$\Big(\log\Big(\frac{1}{r}\Big)\Big)^{-\alpha},$&
			$m_0 = 1$&
			$s \geq 2$
			\\
			$\alpha\in(0,\,\infty)$&
			&
			\\[5pt]
			\bottomrule
		\end{tabulary}
		\caption{Some examples for the index $s$ of $\Zyg^s$ for different moduli of continuity.}
		\label{table04}
	\end{table}

	We also refer to Corollaries~\ref{Intro:cor1} and \ref{Intro:cor2} in Section~\ref{Intro} as examples of our result for $\rho(r) = r$.

	\section{Definitions and tools}\label{DEF}
	Let $x = (x_1,\,\,\ldots,\,x_n)$ be the variables in the $n$-dimensional Euclidean space $\R^n$ and by $\xi = (\xi_1,\,\ldots,\,\xi_n)$ we denote the dual variables. Furthermore, we set $\jbl\xi\jbr^2 = 1 + |\xi|^2$.
	We use the standard multi-index notation. Precisely, let $\Z$ be the set of all integers and $\Z_+$ the set of all non-negative integers. Then $\Z^n_+$ is the set of all $n$-tuples $\alpha = (\alpha_1,\,\ldots,\,\alpha_n)$ with $a_k \in \Z_+$ for each $k = 1,\,\ldots,\,n$. The length of $\alpha \in \Z^n_+$ is given by $|\alpha| = \alpha_1 + \ldots + \alpha_n$.\\
	Let $u = u(t,\,x)$ be a differentiable function, we then write
	\begin{equation*}
	u_t(t,\,x) = \partial_t u (t,\,x) = \frac{\partial}{\partial t} u(t,\,x),
	\end{equation*}
	and
	\begin{equation*}
	\partial_x^\alpha u (t,\,x) = \left(\frac{\partial}{\partial x_1}\right)^{\alpha_1} \ldots\left(\frac{\partial}{\partial x_n}\right)^{\alpha_n} u(t,\,x).
	\end{equation*}
	Using the notation $D_{x_j} = -\I \frac{\partial}{\partial x_j}$, where $\I$ is the imaginary unit, we write also
	\begin{equation*}
	D_x^\alpha = D_{x_1}^{\alpha_1} \cdots D_{x_n}^{\alpha_n}.
	\end{equation*}
	Similarly, for $x\in \R^n$ we set
	\begin{equation*}
	x^\alpha = x_1^{\alpha_1} \cdots x_n^{\alpha_n}.
	\end{equation*}
	
	Let $f$ be a continuous function in an open set $\Omega \subset \R^n$. By $\supp f$ we denote the support of $f$, i.e. the closure in $\Omega$ of $\{x \in \Omega\,|\,f(x) \neq 0\}$. By $C^k(\Omega)$, $0 \leq k \leq \infty$, we denote the set of all functions $f$ defined on $\Omega$, whose derivatives $\partial^\alpha_x f$ exist and are continuous for $|\alpha| \leq k$. By $\Czi(\Omega)$ we denote the set of all functions $f \in C^\infty(\Omega)$ having compact support in $\Omega$. The Sobolev space $H^{k,p}(\Omega)$ consists of all functions that are $k$ times differentiable in Sobolev sense and have (all) derivatives in $L^p(\Omega)$.
	
	We use $C$ as a generic positive constant which may be different even in the same line.
	
	An import tool in our approach is the division of the extended phase space into zones. For this purpose we define $t_\xi$ by
	\begin{equation*}
	t_\xi = N \eta(|\xi|^{-1}),\,N \geq 2,\,|\xi| \geq M.
	\end{equation*}
	The pseudodifferential zone $\Zpd(N,\,M)$ is then given by
	\begin{equation*}
	\Zpd(N,\,M) = \{(t,\,x,\,\xi) \in [0,\,T]\times\R^n_x\times\R^n_\xi\,:\,t \leq t_\xi,\,|\xi| \geq M\},
	\end{equation*}
	consequently, the hyperbolic zone $\Zhyp(N,\,M)$ is defined by
	\begin{equation*}
	\Zhyp(N,\,M) = \{(t,\,x,\,\xi) \in [0,\,T]\times\R^n_x\times\R^n_\xi\,:\,t \geq t_\xi,\,|\xi| \geq M\}.
	\end{equation*}

	Let us recall some results and definitions for the symbol space  $\Zyg^s\Sy^m_{1,\,0}$ and related pseudodifferential operators.
	
	\subsection{Pseudodifferential operators with limited smoothness}\label{sec:pseudo}
	
	We introduce the standard symbol space $\Sy^m_{\rho,\,\delta}$ as in \cite{Hormander.2007} and the space $\Zyg^s\Sy^m_{1,\,\delta}$ of symbols having limited smoothness in $x$.
	
	\begin{Definition}
		Let $m,\,\rho,\,\delta \in \R$ with $0 < \rho \leq 1$ and $0 \leq \delta <1$. A function $p = p (x,\,\xi) \in C^\infty(\R^{2n})$ belongs to $\Sy^m_{\rho,\,\delta}$ if for all $\alpha,\,\beta \in \N^n$ we have the estimate
		\begin{equation*}
		|\partial_\xi\alpha D_x^\beta p(x,\,\xi)| \leq C_{\alpha,\,\beta} \jxi^{m-|\alpha|\rho + |\beta|\delta},
		\end{equation*}
		for all $\xi \in \R^n$.
	\end{Definition}
	
	\begin{Definition}
		Let $s,\,m,\,\delta \in \R$ with $s \geq 0,\,0 \leq \delta < 1$. Then we denote by $\Zyg^s\Sy^m_{1,\,\delta}$ the set of all functions $p = p(x,\,\xi)$ which are smooth in $\xi$ and belong to the Zygmund space $\Zyg^s$ (Definition~\ref{Def:Zyg}) with respect to $x$ such that for all $\alpha,\,\beta \in \N^n$ with $|\beta| < [s]$, we have the estimates
		\begin{equation*}
		|\partial_\xi^\alpha D_x^\beta p(x,\,\xi)| \leq C_{\alpha,\,\beta} \jxi^{m-|\alpha|+|\beta|\delta}, \text{ for } |\beta| \leq [s],
		\end{equation*}
		and
		\begin{equation*}
		\|\partial_\xi^\alpha p(\cdot,\,\xi)\|_{\Zyg^s} \leq C_{\alpha,\,s} \jxi^{m-|\alpha|+\delta s},
		\end{equation*}
		where $[s]$ denotes the largest integer strictly smaller than $s$.
	\end{Definition}

	\subsubsection{Mapping properties}
	The following mapping results as well as the results for composition, adjoint and sharp G{\aa}rding's inequality for pseudodifferential operators with limited smoothness can be proved using a technique called symbol smoothing. We refer the reader to Section~1.3 in \cite{Taylor.1991} for information about that technique.
	
	Let us now briefly review some mapping properties of operators from $\OPS^m_{1,\,\delta}$ and $\Zyg^s\OPS^m_{1,\,\delta}$.
	
	\begin{Proposition}[Chapter 3, Theorem 2.7. in \cite{Kumanogo.1982}]
		Let $p(x,\,D_x) \in \OPS^m_{\rho,\,\delta}$ with $0 \leq \delta < \rho \leq 1$, then
		\begin{equation*}
		p(x,\,D_x): H^{r+m} \rightarrow H^r,
		\end{equation*}
		continuously, for all $r \in \R$.	
	\end{Proposition}

	The following results are valid for the more general nonhomogeneous Besov spaces $B^s_{p,\,q}$.
	
	\begin{Definition}[\cite{Bahouri.2011}]
		Let $s \in \R,\,1\leq p,q \leq \infty$. The nonhomogeneous Besov space $B^s_{p,\,q}$ consists of all tempered distributions $u$ such that
		\begin{equation*}
		\| u\|_{B^s_{p,\,q}} := \Big\| \Big(2^{js} \|\Delta_j u\|_{L^p}\Big)_{j\in \Z}\Big\|_{l^q(\Z)} < \infty.
		\end{equation*}
	\end{Definition}	
	
	\begin{Corollary}
		We note that for $s > 0$ we have $B^s_{\infty,\,\infty} = \Zyg^s$ and for $s \in \R$ we have $B^s_{2,\,2} = H^s$.
	\end{Corollary}
	
	\begin{Proposition}[Lemma 3.4 in \cite{Abels.2005}]\label{Prop:Map2}
		Let $s > 0,\, 0 \leq \delta \leq 1,\, 1 \leq l,\,q \leq \infty,\,m\in\R$, and $p(x,\,D_x) \in \Zyg^s\OPS^m_{1,\,\delta}$, then
		\begin{align*}
		p(x,\,D_x)&: B^{r+m}_{q,\,l} \rightarrow B^r_{q,\,l},\\
		\end{align*}
		if $-(1-\delta)s < r < s$.
	\end{Proposition}
	
	\begin{Proposition}[Lemma 3.5 in \cite{Abels.2005}]\label{Prop:Map1}
		Let $s_1,\,s_2 > 0,\, 1 \leq l,\,q \leq \infty,\,m\in\R$. Let $p(x,\,D_x) \in \Zyg^{s_1}\OPS^m_{1,\,0}\cap \Zyg^{s_2}\OPS^{m-\vartheta}_{1,\,0}$ for some $\vartheta \in (0,\,s_1)$, then
		\begin{align*}
		p(x,\,D_x)&: B^{r+m-\vartheta}_{q,\,l} \rightarrow B^r_{q,\,l},\\
		\end{align*}
		if $-s_1+\vartheta < r < s_1$.
	\end{Proposition}
	
	\subsubsection{Composition, adjoint and sharp G{\aa}rding's inequality}\label{sec:pseudo:calc}
	
	For two symbols $p_1$ and $p_2$ we introduce the notation
	\begin{equation*}
	(p_1\#_k p_2)(x,\,\xi) := \sum\limits_{|\alpha| \leq k} \frac{1}{\alpha!} \partial_\xi^\alpha p_1(x,\,\xi) D_x^\alpha p_2(x,\,\xi),
	\end{equation*}
	for $k \in \N$. Consequently, we write $(p_1\#_k p_2)(x,\,D_x) =\OP(p_1 \#_k p_2)(x,\,D_x)$.

	\begin{Proposition}[Theorem 3.6 in \cite{Abels.2005}]\label{Prop:ZygComp}
		Let $1\leq p,\,q \leq \infty,\,m_1,\,m_2 \in \R,\,s_1,\,s_2 > 0$, choose  $\vartheta \in (0,\,s_2)$ and set $s = \min\{s_1,\,s_2 - [\vartheta]\}$. Let $p_1=p_1(x,\,D_x) \in \Zyg^{s_1}\OPS^{m_1}_{1,\,0}$ and $p_2=p_2(x,\,D_x) \in \Zyg^{s_2}\OPS^{m_2}_{1,\,0}$. For every $r$ such that
		\begin{equation*}
		|r| < s, \qquad r > -(s_2-\vartheta), \qquad -s_2+\vartheta < r + m_1 < s_2,
		\end{equation*}
		we have that
		\begin{equation*}
		R_\vartheta = R_\vartheta(x,\,D_x) := p_1(x,\,D_x) p_2(x,\,D_x) - (p_1 \#_{[\vartheta]}p_2)(x,\,D_x),
		\end{equation*}
		is a bounded operator from
		\begin{equation*}
		B^{r+m_1+m_2 - \vartheta}_{p,\,q} \rightarrow B^r_{p,\,q}.
		\end{equation*}
		The analogous result holds for Bessel potential spaces if $1 <p< \infty$ and $\vartheta \notin \N$.
	\end{Proposition}
	\begin{Remark}
		We cite \cite{Abels.2005} not because it is the first result of this kind but because the technique and notations used there are similar to ours. An earlier version of this result can also be found in \cite{Marschall.1987}.
	\end{Remark}
	
	Let us briefly consider the cases of Proposition~\ref{Prop:ZygComp} that are vital for our approach.	
	\begin{Corollary}\label{Cor:ZygComp}
		Let $p = q = 2$ and take
		\begin{itemize}
			\item $p_1,\,p_2 \in \Zyg^{1+\ve}\OPS^0_{1,\,0}$, then
			\begin{equation*}
			R_{\vartheta=1}: H^{r-1} \rightarrow H^r, \text{ if } -\ve < r < \ve;
			\end{equation*}
			\item $p_1 \in \Zyg^{1+\ve}\OPS^1_{1,\,0},\,p_2 \in \Zyg^{1+\ve}\OPS^0_{1,\,0}$, then
			\begin{equation*}
			R_{\vartheta=1}: H^{r} \rightarrow H^r, \text{ if } -\ve < r < \ve;
			\end{equation*}
			\item $p_1 \in \Zyg^{1+\ve}\OPS^0_{1,\,0},\,p_2 \in \Zyg^{1+\ve}\OPS^1_{1,\,0}$, then
			\begin{equation*}
			R_{\vartheta=1}: H^{r} \rightarrow H^r,  \text{ if } -\ve < r < \ve.
			\end{equation*}
		\end{itemize}
	\end{Corollary}
	
	The following proposition states a result about the adjoint of an operator from $\Zyg^s\OPS^1_{1,\,0}$.
	
	\begin{Proposition}[Proposition 2.3.A in \cite{Taylor.1991}]\label{Prop:ZygAdj}
		Let $p = p(x,\,D_x) \in \Zyg^s\OPS^1_{1,\,0}$ with $s >1$. Then we have that
		\begin{equation*}
		R = R(x,\,D_x) := p(x,\,D_x)^\ast - q(x,\,D_x),
		\end{equation*}
		is a bounded operator from
		\begin{equation*}
		H^r \rightarrow H^r,
		\end{equation*}
		if $1-s < r < s$, where
		\begin{equation*}
		\sigma(q(x,\,D_x))(x,\,\xi) = \xbar{p(x,\,\xi)}.
		\end{equation*}
	\end{Proposition}
	\begin{Remark}
		In his book \cite{Taylor.1991} Taylor has the condition $-s < r < s$. However, I am not sure that this is correct. From the short comments he gives about the proof, I can only obtain $1-s < r < s$.
	\end{Remark}	
	
	Lastly, the following proposition gives a sharp G{\aa}rding inequality for operators with symbols in $\Zyg^s\Sy^m_{1,\,0}$.
	\begin{Proposition}[Corollary II.5 in \cite{Tataru.2002}]\label{Prop:ZygGaar}
		Consider the $N\times N$ symbol $p(x,\,\xi) \in \Zyg^s\Sy^m_{1,\,0}$  with $p(x,\,\xi) \geq 0$. Then for all $u \in \Czi$, we have
		\begin{equation*}
		\Re(p(x,\,D_x)u,\,u) \geq - C_1 \|u\|_{L^2}^2,
		\end{equation*}
		provided that $0 \leq s \leq 2$ and $m \leq \frac{2s}{s+2}$.
	\end{Proposition}
	\begin{Remark}\label{Rem:ZygGaar}
		If we have $p(x,\,\xi) \in \Zyg^s\Sy^1_{1,\,0}$, then the condition
		\begin{equation*}
		1\leq \frac{2s}{s+2},
		\end{equation*}
		yields $s \geq 2$.
	\end{Remark}

	\section{Proof}\label{PROOF}
	
	The steps of the proof are basically the same as in \cite{Kinoshita.2005}. We just have to pay attention to the fact that the used pseudodifferential operators and symbols are not smooth with respect to $x$.
	
	In Section~\ref{sec:Reg} we introduce regularized coefficients and characteristic roots which are smooth with respect to time. After deriving some estimates for the regularized roots, we introduce a suitable symbol class in Section~\ref{sec:sym} which takes account of the behavior of the regularized roots in each zone of the extended phase space. We continue by transforming the original Cauchy problem to a Cauchy problem for a first order system in Section~\ref{sec:trans} and perform two steps of diagonalization in Section~\ref{sec:dia}. Finally, after another change of variables to deal with some lower order terms in Section~\ref{sec:con}, we conclude the proof in Section~\ref{sec:conclusion}.

	As written above, we divide the extended phase space into two zones. The basic idea is that in the pseudodifferential zone $\Zpd(N,\,M)$ we use the global condition on the coefficients, whereas in the hyperbolic zone $\Zhyp(N,\,M)$ we use the local conditions on the coefficients. The separating line of these zones is given by
	\begin{equation}\label{Def:zones}
		t_\xi = N \eta(|\xi|^{-1}),\,N \geq 2,\,|\xi| \geq M.
	\end{equation}
	The pseudodifferential zone $\Zpd(N,\,M)$ is then given by
	\begin{equation*}
		\Zpd(N,\,M) = \{(t,\,x,\,\xi) \in [0,\,T]\times\R^n_x\times\R^n_\xi\,:\,t \leq t_\xi,\,|\xi| \geq M\},
	\end{equation*}
	consequently, the hyperbolic zone $\Zhyp(N,\,M)$ is defined by
	\begin{equation*}
		\Zhyp(N,\,M) = \{(t,\,x,\,\xi) \in [0,\,T]\times\R^n_x\times\R^n_\xi\,:\,t \geq t_\xi,\,|\xi| \geq M\}.
	\end{equation*}

	\subsection{Regularization} \label{sec:Reg}
	Since the coefficients are not smooth with respect to time in $t = 0$, it is helpful to regularize them.
	\begin{Definition}\label{Def:Regularization}
		Let $\psi \in \Czi(\R)$ be a given function satisfying $\int_\R \psi(x) \rmd x = 1$ and $0 \leq \psi(x) \leq 1$ for any $x \in \R$  with $\supp \psi \subset [{-1},1]$. Let $\ve > 0$ and set $\psi_\ve(x) = \frac{1}{\ve} \psi\left(\frac{x}{\ve}\right)$.
		Then we define
		\begin{equation*}
			a_{\ve,\,m-j,\,\gamma}(t,\,x) := (a_{m-j,\,\gamma}\ast_t \psi_\ve)(t,\,x),
		\end{equation*}
		for $j = 0,\,\ldots,\,m-1$.
	\end{Definition}
	For convenience, we write
	\begin{equation*}
		a(t,\,x,\,\xi) := \sum\limits_{j=0}^{m-1}\sum\limits_{|\gamma|+j = m} a_{m-j,\,\gamma}(t,\,x)\frac{\xi^\gamma}{|\xi|^{m-j}},
	\end{equation*}
	and
	\begin{equation*}
		a_\ve(t,\,x,\,\xi) :=  \sum\limits_{j=0}^{m-1}\sum\limits_{|\gamma|+j = m} a_{\ve,\,m-j,\,\gamma}(t,\,x)\frac{\xi^\gamma}{|\xi|^{m-j}}.
	\end{equation*}
	\begin{Proposition}[\cite{Kinoshita.2005,Cicognani.2017}]\label{Prop:RegCoeff}
		We choose $\ve = \jxi^{-1}$. Then the regularized coefficients satisfy the following estimates for all $\alpha \in \N^n$.
		\begin{enumerate}[label = (\roman*),align = left, leftmargin=*]
			\item For $(t,\,\xi) \in [0,\,T]\times\{|\xi|\geq M\}$, we have \begin{equation*}
				\Big\| \partial_\xi^\alpha a_\ve(t,\,\cdot,\,\xi)\Big\|_{\Zyg^s} \leq C_{\alpha,\,s} \jxi^{-|\alpha|}.
			\end{equation*}
			\item For $(t,\,\xi) \in [0,\,T]\times\{|\xi|\geq M\}$, we have \begin{equation*}
				\Big\| \partial_\xi^\alpha\Big( a_\ve(t,\,\cdot,\,\xi)-a(t,\,\cdot,\,\xi)\Big)\Big\|_{\Zyg^s} \leq C_{\alpha,\,s} \jxi^{-1-|\alpha|} \frac{1}{\eta(|\xi|^{-1})}.
			\end{equation*}
			\item For $(t,\,\xi) \in [t_\xi,\,T]\times\{|\xi|\geq M\}$, we have \begin{equation*}
				\Big\| \partial_\xi^\alpha\Big( a_\ve(t,\,\cdot,\,\xi)-a(t,\,\cdot,\,\xi)\Big)\Big\|_{\Zyg^s} \leq -C_{\alpha,\,s} \jxi^{-1-|\alpha|} \varrho(|\xi|^{-1}) \frac{\rmd}{\rmd t}\Big(\frac{1}{\varrho(\eta^{-1}(t-|\xi|^{-1}))}\Big).
			\end{equation*}
			\item For $(t,\,\xi) \in [0,\,T]\times\{|\xi|\geq M\}$, we have
			\begin{equation*}
				\Big\| \partial_\xi^\alpha\partial_t a_\ve(t,\,\cdot,\,\xi)\Big\|_{\Zyg^s} \leq C_{\alpha,\,s} \jxi^{-|\alpha|}  \frac{1}{\eta(|\xi|^{-1})}.
			\end{equation*}
			\item For $(t,\,\xi) \in [t_\xi,\,T]\times\{|\xi|\geq M\}$, we have
			\begin{equation*}
				\Big\|\partial_\xi^\alpha\partial_t a_\ve(t,\,\cdot,\,\xi)\Big\|_{\Zyg^s} \leq C_{\alpha,\,s} \jxi^{-|\alpha|}  \Big(-\frac{\rmd}{\rmd t}\Big(\frac{1}{\eta^{-1}(t-|\xi|^{-1})}\Big)\Big)^{\frac{1}{2}}.
			\end{equation*}
			\item For $(t,\,\xi) \in [t_\xi,\,T]\times\{|\xi|\geq M\}$, we have \begin{equation*}
				\Big\| \partial_\xi^\alpha \partial_t^2 a_\ve(t,\,\cdot,\,\xi)\Big\|_{\Zyg^s} \leq -C_{\alpha,\,s} \jxi^{1-|\alpha|} \varrho(|\xi|^{-1}) \frac{\rmd}{\rmd t}\Big(\frac{1}{\varrho(\eta^{-1}(t-|\xi|^{-1}))}\Big).
			\end{equation*}
		\end{enumerate}
	\end{Proposition}
	\begin{Proof}
		The estimates for $(t,\,\xi) \in [0,\,T]\times\{|\xi| \geq M\}$ follow from Proposition 4.3 in \cite{Cicognani.2017}. The estimates for $(t,\,\xi) \in [t_\xi]\times\{|\xi|\geq M\}$ are derived following the proofs of Lemma 4.1 and Lemma 5.1. in \cite{Kinoshita.2005}.
	\end{Proof}
	
	We introduce $\lambda_k = \lambda_k(t,\,x,\,\xi)$ to be the solutions to
	\begin{equation*}
		\lambda^m - \sum\limits_{j = 0}^{m-1} \sum\limits_{|\gamma|+j = m} a_{\ve,\,m-j,\,\gamma}(t,\,x) \xi^\gamma \lambda^j = 0.
	\end{equation*}
	We renumber these regularized roots such that $\lambda_1 < \ldots < \lambda_m$.
	
	\begin{Proposition}[\cite{Kinoshita.2005,Cicognani.2017}]\label{Prop:RegRoots} We choose $\ve = \jxi^{-1}$. Then the regularized roots satisfy the following relations for all $\alpha \in \N^n$ and all $k = 1,\,\ldots,\,m$.
		\begin{enumerate}[label = (\roman*),align = left, leftmargin=*]
			\item We have $\lambda_k \in C([0,\,T];\,\Zyg^s\Sy^1_{1,\,0})$.
			\item For $(t,\,\xi) \in [0,\,T]\times\{|\xi|\geq M\}$, we have \begin{equation*}
				\Big\| \partial_\xi^\alpha\Big( \lambda_k(t,\,\cdot,\,\xi)-\tau_k(t,\,\cdot,\,\xi)\Big)\Big\|_{\Zyg^s} \leq C_{\alpha,\,s} \jxi^{-|\alpha|} \frac{1}{\eta(|\xi|^{-1})}.
			\end{equation*}
			\item For $(t,\,\xi) \in [t_\xi,\,T]\times\{|\xi|\geq M\}$, we have \begin{equation*}
				\Big\| \partial_\xi^\alpha\Big( \lambda_k(t,\,\cdot,\,\xi)-\tau_k(t,\,\cdot,\,\xi)\Big)\Big\|_{\Zyg^s} \leq -C_{\alpha,\,s} \jxi^{-|\alpha|} \varrho(|\xi|^{-1}) \frac{\rmd}{\rmd t}\Big(\frac{1}{\varrho(\eta^{-1}(t-|\xi|^{-1}))}\Big).
			\end{equation*}
			\item For $(t,\,\xi) \in [0,\,T]\times\{|\xi|\geq M\}$, we have
			\begin{equation*}
				\Big\| \partial_\xi^\alpha\partial_t \lambda_k(t,\,\cdot,\,\xi)\Big\|_{\Zyg^s} \leq C_{\alpha,\,s} \jxi^{1-|\alpha|}  \frac{1}{\eta(|\xi|^{-1})}.
			\end{equation*}
			\item For $(t,\,\xi) \in [t_\xi,\,T]\times\{|\xi|\geq M\}$, we have
			\begin{equation*}
				\Big\|\partial_\xi^\alpha\partial_t \lambda_k(t,\,\cdot,\,\xi)\Big\|_{\Zyg^s} \leq C_{\alpha,\,s} \jxi^{1-|\alpha|}  \Big(-\frac{\rmd}{\rmd t}\Big(\frac{1}{\eta^{-1}(t-|\xi|^{-1})}\Big)\Big)^{\frac{1}{2}}.
			\end{equation*}
			\item For $(t,\,\xi) \in [t_\xi,\,T]\times\{|\xi|\geq M\}$, we have \begin{equation*}
				\Big\| \partial_\xi^\alpha \partial_t^2 \lambda_k(t,\,\cdot,\,\xi)\Big\|_{\Zyg^s} \leq -C_{\alpha,\,s} \jxi^{2-|\alpha|} \varrho(|\xi|^{-1}) \frac{\rmd}{\rmd t}\Big(\frac{1}{\varrho(\eta^{-1}(t-|\xi|^{-1}))}\Big).
			\end{equation*}
		\end{enumerate}
	\end{Proposition}
	
	\subsection{Symbol space} \label{sec:sym}
	We introduce the symbol space $P^m$ which takes account of the estimates we observe in Proposition~\ref{Prop:RegRoots}.
	\begin{Definition}
		Let $b = b(t,\,x,\,\xi) \in L^\infty([0,\,T];\,\Zyg^s\Sy^{m+1}_{1,\,0})$. Then $b$ belongs to the symbol class $P^m$, if it satisfies
		\begin{equation*}
			\Big\|\partial_\xi^\alpha b(t,\,\cdot,\,\xi)\Big\|_{\Zyg^s} \leq C_{\alpha,\,s} |\xi|^{m-|\alpha|} \frac{1}{\eta(|\xi|^{-1})},
		\end{equation*}
		for $(t,\,\xi) \in [0,\,t_\xi]\times\{|\xi| \geq M\}$,	and
		\begin{equation*}
			\begin{aligned}
				\Big\|\partial_\xi^\alpha b(t,\,\cdot,\,\xi)\Big\|_{\Zyg^s} \leq C_{\alpha,\,s} \jxi^{m-|\alpha|} \varrho(|\xi|^{-1}) \Big(-\frac{\rmd}{\rmd t}\Big(\frac{1}{\varrho(\eta^{-1}(t-|\xi|^{-1}))}\Big)\Big)\\
				+ C_{\alpha,\,s} \jxi^{m-|\alpha|-1}  \Big(-\frac{\rmd}{\rmd t}\Big(\frac{1}{\eta^{-1}(t-|\xi|^{-1})}\Big)\Big),
			\end{aligned}
		\end{equation*}
		for $(t,\,\xi) \in [t_\xi,\,T]\times\{|\xi| \geq M\}$.
	\end{Definition}
	
	\begin{Corollary}[Lemma 4.3 in \cite{Kinoshita.2005}]
		If $b \in P^m$, then $b \in L^\infty([0,\,T];\,\Zyg^s\Sy^{m+1}_{1,\,0})$. Furthermore, if $\psi \in \Czi(\R^n)$ with $\supp \psi \subset \{|\xi|\geq M\}$, then \begin{equation*}
			\psi(\xi) \lambda_k(t,\,x,\,\xi) \in L^\infty([0,\,T];\,\Zyg^s\Sy^1_{1,\,0}),
		\end{equation*}
		and 
		\begin{equation*}
			\psi(\xi)D_t \lambda_k(t,\,x,\,\xi) \in L^\infty([0,\,T];\,\Zyg^s\Sy^2_{1,\,0}).
		\end{equation*}
		
	\end{Corollary}
	
	With these preparations complete, let us now tend to the original Cauchy problem.

	\subsection{Transformation to a first-order system} \label{sec:trans}
	
	We consider the Cauchy problem \eqref{CP:eq1} and set
	\begin{equation*}
		U = (\jbl D_x \jbr^{m-1} \psi(D_x) u,\,\jbl D_x \jbr^{m-2} \psi(D_x)D_t u,\,\ldots,\,\psi(D_x)D_t^{m-1}u)^T,
	\end{equation*}
	where $\psi$ localizes to large frequencies.
	This yields
	\begin{equation*}
		D_t U = A(t,\,x,\,D_x) U + R(t,\,x,\,D_x) U,
	\end{equation*}
	with
	\begin{equation*}
		A(t,\,x,\,\xi) = \sigma(A) = \begin{pmatrix}
			0 & \jxi& 0&\cdots&0\\
			\vdots& & \ddots& &\vdots\\
			0& \cdots& 0& \jxi&0\\
			a_{m,\,\gamma}(t,\,x,\,\xi) & \cdots&a_{m-j,\,\gamma}(t,\,x,\,\xi) & \cdots & a_{1,\,\gamma}(t,\,x,\,\xi)
		\end{pmatrix},
	\end{equation*}
	where
	\begin{equation}\label{amj}
	a_{m-j,\,\gamma}(t,\,x,\,\xi) = \sum\limits_{|\gamma| = m-j} a_{m-j,\,\gamma}(t,\,x) \xi^\gamma \jxi^{-(m-1-j)} \in L^\infty([0,\,T];\,\Zyg^s\Sy^1_{1,\,0}),
	\end{equation} and $R(t,\,x,\,D_x) \in L^\infty([0,\,T];\,\OPS^{-\infty}_{1,\,0})$.
	
	\subsection{Diagonalization}\label{sec:dia}
	
	We perform two steps of diagonalization. The first step is done in both zones, whereas the second step of diagonalization is only done in the hyperbolic zone.
	
	\subsubsection{First step of diagonalization}
	
	We introduce the diagonalizer $M_1 = M_1(t,\,x,\,D_x)$ with symbol
	\begin{equation*}
		M_1(t,\,x,\,\xi) = \sigma (M_1) = \begin{pmatrix}
			1 & \cdots &1\\
			\frac{\lambda_1}{\jxi} & \cdots & \frac{\lambda_m}{\jxi}\\
			\vdots & & \vdots\\
			\Big(\frac{\lambda_1}{\jxi}\Big)^{m-1} & \cdots & \Big(\frac{\lambda_m}{\jxi}\Big)^{m-1}
		\end{pmatrix},
	\end{equation*}
	as well as the matrix pseudodifferential operator $\widetilde M_1 = \widetilde M_1(t,\,x,\,D_x)$ with its symbol
	\begin{equation*}
		\sigma(\widetilde M_1) =  \sigma(M_1)^{-1} = (c_{p,\,q}(t,\,x,\,\xi))_{1\leq p ,\,q \leq m},
	\end{equation*}
	given by
	\begin{equation*}
		c_{p,\,q} = (-1)^{q-1} \jxi^{q-1} \sum\limits_{S^{(m)}_{\{p\}}(m-q)} \lambda_{i_1} \cdot \ldots \cdot \lambda_{i_{m-q}} \prod\limits_{\substack{i=1\\i\neq p}}^m (\lambda_i -\lambda_p)^{-1},
	\end{equation*}
	for $1 \leq q \leq m-1$ and by
	\begin{equation*}
		c_{p,\,m} = (-1)^{m-1} \jxi^{m-1} \prod\limits_{\substack{i=1\\i\neq p}}^m (\lambda_i -\lambda_p)^{-1},
	\end{equation*}
	where
	\begin{equation*}
		S^{(m)}_{B}(k) := \{(i_1,\,\ldots,\,i_k) \in \N^k\,:\, 1 \leq i_1 < \ldots < i_k \leq m,\,i_l \notin B,\,l =1 ,\,\ldots,\,k\}.
	\end{equation*}
	
	\begin{Proposition}\label{Dia:M1}
		For the matrix pseudodifferential operators $M_1$ and $\widetilde M_1$ defined above, we have the following properties.
		\begin{enumerate}[label = (\roman*)]
			\item We have $\sigma(M_1),\,\sigma(\widetilde M_1) \in L^\infty([0,\,T];\,\Zyg^s\Sy^0_{1,\,0})$ for $|\xi| \geq M$.
			\item In $\Zhyp(N,\,M) \cup \Zpd(N,\,M)$, we have for $s = 1+ \ve > 1$.
			\begin{equation*}
				\widetilde M_1(t,\,x,\,D_x) M_1(t,\,x,\,D_x) = I + R_1(t,\,x,\,D_x),
			\end{equation*} with $R_1(t,\,x,\,D_x) : H^{r-1} \rightarrow H^r$ continuously for $|r| < \ve$ and $(t,\,x) \in [0,\,T]\times\R^n$.
		\end{enumerate}
	\end{Proposition}
	\begin{Proof}
		Both assertions follow essentially from the composition result for operators from $\Zyg^s\OPS^m_{1,\,0}$ (see Proposition~\ref{Prop:ZygComp} and Corollary~\ref{Cor:ZygComp}) and the observation that $\sigma(\widetilde M_1)$ is the inverse matrix of $\sigma(M_1)$.
	\end{Proof}
	
	We set $U = M_1(t,\,x,\,D_x) U_1$ and obtain that
	\begin{equation*}
		D_t U = (D_t M_1) U_1 + M_1 D_t U_1 = A M_1 U_1 + R U,
	\end{equation*}
	which gives
	\begin{equation}\label{CP:eq2}
	D_t U_1 = \widetilde M_1 A M_1 U_1 - \widetilde M_1 (D_t M_1) U_1 - R_1 (D_t U_1) + \widetilde M_1 R M_1 U_1,
	\end{equation}
	where $R \in L^\infty([0,\,T]; \OPS^{-\infty})$ and $R_1(t,\,x,\,D_x) : H^{r-1} \rightarrow H^r$ continuously for $|r| < \ve$ and $(t,\,x) \in [0,\,T]\times\R^n$.
	
	\begin{Proposition}
		The Cauchy problem \eqref{CP:eq2} is for $s = 1+ \ve > 1$ equivalent to
		\begin{equation}\label{CP:eq3}
		D_t U_1 = (A_1 + B_1 - C_1 + R_2+R_\infty)U_1 - R_1 (D_t U_1),
		\end{equation}
		where
		\begin{equation*}
			\sigma(A_1) = \begin{pmatrix}
				\lambda_1 & &\\
				& \ddots & \\
				&& \lambda_m
			\end{pmatrix}, \qquad \sigma(B_1) \in P^0,
		\end{equation*}
		and $\sigma(C_1 ) = (e_{p,\,q})_{1\leq p,\,q\leq m}$ with
		\begin{equation*}
			e_{p,\,q} = \begin{cases}
				- D_t \lambda_p \sum\limits_{\substack{i = 1\\i \neq p}}^m \frac{1}{\lambda_i - \lambda_p}, & p =q,\\[2em]
				- D_t \lambda_q \frac{\prod\limits_{\substack{i = 1\\i \neq p,\,q}}^m (\lambda_i - \lambda_q)}{\prod\limits_{\substack{i = 1\\i \neq p}}^m(\lambda_i-\lambda_p)} , & p \neq q,
			\end{cases}
		\end{equation*}
		and $R_1$ as above, $R_2(t,\,x,\,D_x): H^r \rightarrow H^r$ for $|r| < \ve$ and $R_\infty=\widetilde M_1 R M_1 : H^r\rightarrow H^q$ for $|q|,\,|r| < s$ and all $(t,\,x) \in [0,\,T]\times\R^n$.
	\end{Proposition}
	\begin{Proof}
		To obtain \eqref{CP:eq3} from \eqref{CP:eq2}, we compute $\widetilde M_1 A M_1$ and $\widetilde M_1 (D_t M_1)$ and analyze the respective symbols.
		
		Following our composition result for operators from $\Zyg^s\OPS^m_{1,\,0}$ (see Proposition~\ref{Prop:ZygComp} and Corollary~\ref{Cor:ZygComp}), we obtain that
		\begin{equation*}
			\sigma(\widetilde M_1 A M_1) = \sigma(\widetilde M_1) \sigma(A) \sigma(M_1) + \widetilde R,
		\end{equation*}
		with $\widetilde R : H^r \rightarrow H^r$ for $|r| < \ve$. We write
		\begin{equation*}
			\begin{aligned}
				\sigma(A) = \begin{pmatrix}
					0 & \jxi& 0&\cdots&0\\
					\vdots& & \ddots& &\vdots\\
					0& \cdots& 0& \jxi&0\\
					a_{\ve,\,m,\,\gamma} & \cdots&a_{\ve,\,m-j,\,\gamma} & \cdots & a_{\ve,\,1,\,\gamma}
				\end{pmatrix}\\ + \begin{pmatrix}
					&&&&\\
					&&&&\\
					&&&&\\
					(a_{m,\,\gamma}-a_{\ve,\,m,\,\gamma})& \cdots&(a_{m-j,\,\gamma}-a_{\ve,\,m-j,\,\gamma}) & \cdots & (a_{1,\,\gamma}-a_{\ve,\,1,\,\gamma}),
				\end{pmatrix},
			\end{aligned}
		\end{equation*}
		where
		\begin{equation*}
			a_{\ve,\,m-j,\,\gamma}=a_{\ve,\,m-j,\,\gamma}(t,\,x,\,\xi) = \sum\limits_{|\gamma| = m-j} a_{\ve,\,m-j,\,\gamma}(t,\,x) \xi^\gamma \jxi^{-(m-1-j)},
		\end{equation*}
		and $a_{m-j,\,\gamma}$ are given by \eqref{amj}. Noting that
		\begin{equation*}
			(a_{m-j,\,\gamma}-a_{\ve,\,m-j,\,\gamma}) \in P^0, \text{ (due to Proposition~\ref{Prop:RegCoeff})}
		\end{equation*}
		we have
		\begin{equation*}
			\sigma(\widetilde M_1) \sigma(A) \sigma(M_1) = \begin{pmatrix}
				\lambda_1 & &\\
				& \ddots & \\
				&& \lambda_m
			\end{pmatrix} + B_1,
		\end{equation*} 
		with $\sigma(B_1) \in P^0$.
		
		Concerning $\widetilde M_1 (D_t M_1)$, we first look at $(D_t M_1)$.
		Due to Proposition~\ref{Prop:RegRoots}, we have for $(t,\,\xi) \in [0,\,T]\times\{|\xi|\geq M\}$ that
		\begin{equation}\label{roots:est1}
		\Big\| \partial_\xi^\alpha\partial_t \lambda_k(t,\,\cdot,\,\xi)\jxi^{-1}\Big\|_{\Zyg^s} \leq C_{\alpha,\,s} \jxi^{-|\alpha|}  \frac{1}{\eta(|\xi|^{-1})},
		\end{equation}
		and for $(t,\,\xi) \in [t_\xi,\,T]\times\{|\xi|\geq M\}$, we have
		\begin{equation}\label{roots:est2}
		\Big\|\partial_\xi^\alpha\partial_t \lambda_k(t,\,\cdot,\,\xi)\jxi^{-1}\Big\|_{\Zyg^s} \leq C_{\alpha,\,s} \jxi^{-|\alpha|}  \Big(-\frac{\rmd}{\rmd t}\Big(\frac{1}{\eta^{-1}(t-|\xi|^{-1})}\Big)\Big)^{\frac{1}{2}}.
		\end{equation}
		Thus $\sigma(D_t M_1)$ does not belong to $P^0$ (due to the exponent $\frac{1}{2}$) but it clearly belongs to $L^\infty([0,\,T];\,\Zyg^s\Sy^1_{1,\,0})$. Application of the composition result (see Proposition~\ref{Prop:ZygComp} and Corollary~\ref{Cor:ZygComp}), yields
		\begin{equation*}
			\sigma(\widetilde M_1 (D_t M_1)) = \sigma(\widetilde M_1) \sigma(D_t M_1) + \xbar R,
		\end{equation*}
		where $\xbar R : H^r \rightarrow H^r$ for $|r| < \ve$. Finally, computing $\sigma(\widetilde M_1) \sigma(D_t M_1) = (e_{p,\,q})_{1\leq p,\,q\leq m}$ yields
		\begin{equation*}
			e_{p,\,q} = \begin{cases}
				- D_t \lambda_p \sum\limits_{\substack{i = 1\\i \neq p}}^m \frac{1}{\lambda_i - \lambda_p}, & p =q,\\[2em]
				- D_t \lambda_q \frac{\prod\limits_{\substack{i = 1\\i \neq p,\,q}}^m (\lambda_i - \lambda_q)}{\prod\limits_{\substack{i = 1\\i \neq p}}^m(\lambda_i-\lambda_p)} , & p \neq q,
			\end{cases}
		\end{equation*}
		and writing $R_2 = \xbar R + \widetilde R$ concludes this part of the proof.
		
		Finally, setting $R_\infty = \widetilde M_1 R M_1$ and observing that
		\begin{align*}
			M_1&: H^r\rightarrow H^r, \text{ for } |r| < s,\\
			R&: H^r\rightarrow H^q, \text{ for any } r,\,q,\\
			\widetilde M_1 &: H^q \rightarrow H^q, \text{ for } |q| < s,
		\end{align*}
		concludes the proof.
		
	\end{Proof}
	
	\begin{Remark}
		Looking at \eqref{CP:eq3}, we note that the operator $A_1$ is diagonal and $\sigma(B_1) \in P^0$. If $\sigma(C_1)$ belonged to $P^0$, there would be no need for a second step of diagonalization and we could continue with the next step of the proof.
		
		We note that the behavior of $\sigma(C_1)$ is characterized by \eqref{roots:est1} and \eqref{roots:est2}. Its behavior in the pseudodifferential zone, i.e. for $t \leq t_\xi$ is fine - there it satisfies the estimate required to belong to $P^0$. However, its behavior in the hyperbolic zone does not fit into $P^0$, which is why we carry out the second step of diagonalization only in the hyperbolic zone.
	\end{Remark}

	\subsubsection{Second step of diagonalization}
	
	We consider the Cauchy problem \eqref{CP:eq3} and want to diagonalize $C_1$. We follow the standard diagonalization procedure and set
	\begin{equation*}
		\sigma(M_2) = (d_{p,\,q})_{1\leq p,\,q \leq m} = I + \chixi \frac{e_{p,\,p}}{\lambda_p-\lambda_q}.
	\end{equation*}
	In other words
	\begin{equation*}
		d_{p,\,q} = \begin{cases}
			1, & p =q,\\
			- \frac{D_t \lambda_q}{\lambda_p-\lambda_q} \prod\limits_{\substack{i = 1\\i \neq p,\,q}}^m (\lambda_i - \lambda_q)\prod\limits_{\substack{i = 1\\i \neq p}}^m(\lambda_i-\lambda_p)^{-1} , & p \neq q,
		\end{cases}
	\end{equation*}
	for $(t,\,x,\,\xi) \in \Zhyp(N,\,M)$. We introduce $\widetilde M_2$ with symbol $\sigma(\widetilde M_2) = \sigma(M_2)^{-1}$, which exists since all columns are linearly independent of each other.
	
	\begin{Proposition}\label{Dia:M2}
		For the matrix pseudodifferential operators $M_2$ and $\widetilde M_2$ defined above, we have the following properties.
		\begin{enumerate}[label = (\roman*)]
			\item We have $\sigma(M_2),\,\sigma(\widetilde M_2) \in L^\infty([0,\,T];\,\Zyg^s\Sy^0_{1,\,0})$ for $|\xi| \geq M$.
			\item In $\Zhyp(N,\,M) \cup \Zpd(N,\,M)$, we have for $s = 1+ \ve > 1$
			\begin{equation*}
				\widetilde M_2(t,\,x,\,D_x) M_2(t,\,x,\,D_x) = I + R_3(t,\,x,\,D_x),
			\end{equation*} with $R_3(t,\,x,\,D_x) : H^{r-1} \rightarrow H^r$ continuously for $|r| < \ve$ and $(t,\,x) \in [0,\,T]\times\R^n$..
		\end{enumerate}
	\end{Proposition}
	\begin{Proof}
		Both assertions follow essentially from the composition result for operators from $\Zyg^s\OPS^m_{1,\,0}$ (see Proposition~\ref{Prop:ZygComp} and Corollary~\ref{Cor:ZygComp}).
	\end{Proof}
	
	We perform the second step of diagonalization by setting $U_1 = M_2 U_2$ and obtain
	\begin{equation}\label{CP:eq4}
	D_t U_1 = (D_t M_2) U_2 + M_2 D_t U_2 =  (A_1 + B_1 - C_1 + R_2+R_\infty)M_2U_2 - R_1 (D_t U_1).
	\end{equation}
	
	\begin{Proposition}
		The Cauchy problem \eqref{CP:eq4} is for $s = 1+ \ve > 1$ equivalent to
		\begin{equation}\label{CP:eq5}
		D_t U_2 = (A_1 + A_2 + B_2 + R_2+R_\infty) U_2 - R_3 (D_t U_2) - R_1(D_t U_1),
		\end{equation}
		where
		\begin{equation*}
			\sigma(A_1) = \begin{pmatrix}
				\lambda_1 & &\\
				& \ddots & \\
				&& \lambda_m
			\end{pmatrix}, \qquad \sigma(A_2 ) = \begin{pmatrix}
				D_t \lambda_1 \sum\limits_{\substack{i = 1\\i \neq 1}}^m \frac{1}{\lambda_i - \lambda_1} & &\\
				& \ddots & \\
				&& D_t \lambda_m \sum\limits_{\substack{i = 1\\i \neq m}}^m \frac{1}{\lambda_i - \lambda_m}
			\end{pmatrix},
		\end{equation*}
		with $\sigma(B_2) \in P^0$ and
		\begin{itemize}
			\item $R_1(t,\,x,\,D_x) : H^{r-1} \rightarrow H^r$ for $|r| < \ve$,
			\item $R_2(t,\,x,\,D_x) : H^{r} \rightarrow H^r$ for $|r| < \ve$,
			\item $R_3(t,\,x,\,D_x) : H^{r-1} \rightarrow H^r$ for $|r| < \ve$, and
			\item $R_\infty(t,\,x,\,D_x) : H^r\rightarrow H^q$ for $|q|,\,|r| < s$.
		\end{itemize}
	\end{Proposition}
	\begin{Proof}
		The result follows from the standard diagonalization routine (see e.g. \cite{Jachmann.2010}) and the composition result in Proposition~\ref{Prop:ZygComp} and Corollary~\ref{Cor:ZygComp}.
	\end{Proof}
	
	We carry out one more change of variables to derive the system from which we conclude the desired energy estimate.
	
	We define the matrix pseudodifferential operator $M_3 = M_3(t,\,x,\,D_x)$ with symbol
	\begin{equation*}
		\sigma(M_3) = \begin{pmatrix}
			w_1 & &\\
			& \ddots&\\
			&&w_m
		\end{pmatrix},
	\end{equation*}
	where
	\begin{equation*}
		w_p = w_p(t,\,x,\,\xi) =\exp\Bigg(\int\limits_0^t\frac{D_s \lambda_p(s,\,x,\,\xi)}{\sum\limits_{\substack{i=1\\i\neq p}}^m (\lambda_i - \lambda_p)(s,\,x,\,\xi)}\rmd s\Bigg).
	\end{equation*}
	We introduce the matrix pseudodifferential operator $\widetilde M_3 = \widetilde M_3(t,\,x,\,D_x)$ with symbol
	\begin{equation*}
		\sigma(\widetilde M_3) = \sigma(M_3)^{-1} = \begin{pmatrix}
			w_1^{-1} & &\\
			& \ddots&\\
			&&w_m^{-1}
		\end{pmatrix}.
	\end{equation*}
	
	\begin{Proposition}\label{Dia:M3}
		For the matrix pseudodifferential operators $M_3$ and $\widetilde M_3$ defined above, we have the following properties.
		\begin{enumerate}[label = (\roman*)]
			\item We have $\sigma(M_3),\,\sigma(\widetilde M_3) \in L^\infty([0,\,T];\,\Zyg^s\Sy^0_{1,\,0})$ for $|\xi| \geq M$.
			\item In $\Zhyp(N,\,M) \cup \Zpd(N,\,M)$, we have for $s = 1+ \ve > 1$
			\begin{equation*}
				\widetilde M_3(t,\,x,\,D_x) M_3(t,\,x,\,D_x) = I + R_4(t,\,x,\,D_x),
			\end{equation*} with $R_4(t,\,x,\,D_x) : H^{r-1} \rightarrow H^r$ continuously for $|r| < \ve$ and $(t,\,x) \in [0,\,T]\times\R^n$.
		\end{enumerate}
	\end{Proposition}
	\begin{Proof}
		Both assertions follow from the observation that
		\begin{equation}\label{Dia:M3:est}
			\int\limits_0^t\frac{D_s \lambda_p(s,\,x,\,\xi)}{\sum\limits_{\substack{i=1\\i\neq p}}^m (\lambda_i - \lambda_p)(s,\,x,\,\xi)}\rmd s \in L^\infty([0,\,T];\,\Zyg^s \Sy^0_{1,\,0}),
		\end{equation}
		for all $(t,\,x,\,\xi) \in [0,\,T]\times\R^n\times\{|\xi| \geq M \}$ and from the composition result for operators from $\Zyg^s\OPS^m_{1,\,0}$ (see Proposition~\ref{Prop:ZygComp} and Corollary~\ref{Cor:ZygComp}).
		Relation~\eqref{Dia:M3:est} can be proved by splitting the integral into two integrals. One integral from $0$ to $t_\xi$ and one integral from $t_\xi$ to $T$. Using the properties of $\lambda_k,\,k = 1,\,\ldots,\,m$ from Proposition~\ref{Prop:RegRoots} and the definition of the zones~\eqref{Def:zones} then gives \eqref{Dia:M3:est}.
	\end{Proof}
	
	Setting $U_2 = M_3 U_3$ then yields
	\begin{equation}\label{CP:eq6}
	D_t U_2 = (D_t M_3) U_3 + M_3 D_t U_3 =  (A_1 + A_2 + B_2 + R_2+R_\infty) U_2 - R_3 (D_t U_2) - R_1(D_t U_1).
	\end{equation}
	
	\begin{Proposition}
		The Cauchy problem \eqref{CP:eq6} is for $s = 1+ \ve > 1$ equivalent to
		\begin{equation*}
			D_t U_3 = (A_1+B_2+R_2 + R_\infty)U_3 - R_4(D_t U_3) - R_3 (D_t U_2) - R_1(D_t U_1),
		\end{equation*}
		where
		\begin{equation*}
			\sigma(A_1) = \begin{pmatrix}
				\lambda_1 & &\\
				& \ddots & \\
				&& \lambda_m
			\end{pmatrix}, \qquad\sigma(B_2) \in P^0,
		\end{equation*}
		and
		\begin{itemize}
			\item $R_1(t,\,x,\,D_x) : H^{r-1} \rightarrow H^r$ for $|r| < \ve$,
			\item $R_2(t,\,x,\,D_x) : H^{r} \rightarrow H^r$ for $|r| < \ve$,
			\item $R_3(t,\,x,\,D_x) : H^{r-1} \rightarrow H^r$ for $|r| < \ve$,
			\item $R_4(t,\,x,\,D_x) : H^{r-1} \rightarrow H^r$ for $|r| < \ve$, and
			\item $R_\infty(t,\,x,\,D_x) : H^r\rightarrow H^q$ for $|q|,\,|r| < s$.
		\end{itemize}
	\end{Proposition}
	\begin{Proof}
		The statement of the proposition follows from the composition result for operators from $\Zyg^s\OPS^m_{1,\,0}$ (see Proposition~\ref{Prop:ZygComp} and Corollary~\ref{Cor:ZygComp}) and the fact that
		\begin{equation*}
			\sigma(\widetilde M_3) \sigma(D_t M_3) = \sigma(A_2).
		\end{equation*}
	\end{Proof}
	
	\subsection{Conjugation}\label{sec:con}
		We want to apply Duhamel's principle to
		\begin{equation*}
			\begin{cases}
				D_t U_3 -(A_1+B_2+R_2 + R_\infty)U_3 =- R_4(D_t U_3) - R_3 (D_t U_2) - R_1(D_t U_1),\\
				U_3(0,\,x) = G(x),
			\end{cases}
		\end{equation*}
		where we consider the terms $- R_4(D_t U_3) - R_3 (D_t U_2) - R_1(D_t U_1)$ as an inhomogeneity and where $G(x)$ denotes the vector containing the transformed initial data. In that way, we can argue that it is sufficient to consider the homogeneous problem
		\begin{equation*}
			D_t W - (A_1+ B_2 + R_2 + R_\infty) W = 0,
		\end{equation*}
		with initial data
		\begin{equation*}
			W(s,\,s,\,x) = - R_4(D_t U_3)(s,\,x) - R_3 (D_t U_2)(s,\,x) - R_1(D_t U_1)(s,\,x).
		\end{equation*}
		We use this approach to obtain an $L^2-L^2$ energy estimate for $W$ without loss of derivatives. We note that it is important for the application of Duhamel's principle that
		\begin{itemize}
			\item $R_1(t,\,x,\,D_x) : H^{r-1} \rightarrow H^r$ for $|r| < \ve$,
			\item $R_3(t,\,x,\,D_x) : H^{r-1} \rightarrow H^r$ for $|r| < \ve$, and
			\item $R_4(t,\,x,\,D_x) : H^{r-1} \rightarrow H^r$ for $|r| < \ve$,
		\end{itemize}
		which ensures that the initial data $- R_4(D_t U_3)(s,\,x) - R_3 (D_t U_2)(s,\,x) - R_1(D_t U_1)(s,\,x)$ is in $L^2$ with respect to $x$.

	Let us now consider the homogeneous problem
	\begin{equation*}
		D_t W - (A_1+ B_2 + R_2 + R_\infty) W = 0,
	\end{equation*}
	and define
	\begin{equation*}
		V:= \exp\bigg(- \int\limits_0^t \vartheta(s,\,D_x) \rmd s\bigg)W,
	\end{equation*}
	with
	\begin{equation*}
		\vartheta(t,\,\xi) = K(2+ \vartheta_0(t,\,\xi)),
	\end{equation*}
	where
	\begin{equation*}
		\begin{aligned}
			\vartheta_0(t,\,\xi) = \bigg(1-\chi\bigg(\frac{t}{2N\eta(\jxi^{-1})}\bigg)\bigg) \frac{1}{\eta(\jxi^{-1})} + \chi\bigg(\frac{t}{N\eta(\jxi^{-1})}\bigg)\bigg[-\varrho(\jxi^{-1}) \frac{\rmd}{\rmd t}\bigg(\frac{1}{\varrho(\eta^{-1}(t-\jxi^{-1}))}\bigg)\\ - \jxi^{-1} \frac{\rmd}{\rmd t}\bigg(\frac{1}{\eta^{-1}(t-\jxi^{-1})}\bigg) \bigg]
		\end{aligned}
	\end{equation*}
	
	\begin{Proposition}[Lemmas 5.6 and 5.7. in \cite{Kinoshita.2005}]\label{Prop:Theta}
		We have
		\begin{enumerate}[label = (\roman*)]
			\item $\vartheta_0(t,\,\xi) \in L^\infty([0,\,T];\,\Sy^{m_0}_{1,\,0}) \subset L^\infty([0,\,T];\,\Sy^1_{1,\,0})$, and
			\item $\int\limits_0^t \vartheta_0(s,\,\xi) \rmd s \in L^\infty([0,\,T];\,\Sy^0_{1,\,0})$.
		\end{enumerate}
	\end{Proposition}
	Applying this transformation yields
	\begin{equation*}
		D_t V = \I \vartheta(t,\,D_x) I V + \exp\bigg(- \int\limits_0^t \vartheta(s,\,D_x) \rmd s\bigg) (A_1 + B_2 + R_2 + R_\infty)\exp\bigg(\int\limits_0^t \vartheta(s,\,D_x) \rmd s\bigg)V,
	\end{equation*}
	which gives
	\begin{equation*}
		\begin{aligned}
			\partial_t V + \vartheta(t,\,D_x) I V - \I(A_1 + B_1 + R_2 + R_\infty)V \\- \I \bigg[\exp\bigg(- \int\limits_0^t \vartheta(s,\,D_x) \rmd s\bigg), A_1+B_2+R_2 + R_\infty\bigg] \exp\bigg( \int\limits_0^t \vartheta(s,\,D_x) \rmd s\bigg)V = 0.
		\end{aligned}
	\end{equation*}
	We define
	\begin{align*}
		Q_0 &:= K(1+ \vartheta_0(t,\,D_x))I - \I(A_1 + B_2+ R_2 + R_\infty),\\
		Q_1 &:= KI- \I \bigg[\exp\bigg(- \int\limits_0^t \vartheta(s,\,D_x) \rmd s\bigg), A_1+B_2+R_2 + R_\infty\bigg] \exp\bigg( \int\limits_0^t \vartheta(s,\,D_x) \rmd s\bigg).
	\end{align*}
	
	\begin{Proposition}[Lemma 4.6 in \cite{Kinoshita.2005}]\label{Prop:Q0}
		For $s > 1$, we have that
		\begin{enumerate}[label = (\roman*)]
			\item $\psi(D_x) Q_0 \in L^\infty([0,\,T];\,\Sy^{m_0}_{1,\,0}) \subset L^\infty([0,\,T];\, \Zyg^s\OPS^1_{1,\,0})$,
			\item $\sigma\Big(\psi(D_x) \frac{Q_0 + Q_0^\ast}{2}\Big) \geq \vartheta_0(t,\,\xi)I$.
		\end{enumerate}
	\end{Proposition}
	\begin{Proof}
		The first statement follows from the structure of $Q_0$ and Proposition~\ref{Prop:Theta}. For the second statement we employ Proposition~\ref{Prop:ZygAdj} to deal with adjoints of operators from $\Zyg^s\OPS^1_{1,\,0}$. We obtain that
		\begin{align*}
			\sigma(A_1^\ast) &=  \xbar{\sigma(A_1)} + r_1,\\
			\sigma(B_2^\ast) &=  \xbar{\sigma(B_2)} + r_2,
		\end{align*}
		with $r_1,\,r_2: L^2 \rightarrow L^2$ continuously. Since the problem is strictly hyperbolic, we conclude that
		\begin{equation*}
			\I A_1 = - (\I A_1)^\ast - r_1.
		\end{equation*}
		Furthermore, since $\sigma(B_2) \in P^0$, we know that
		\begin{equation*}
			|\sigma(B_2)| \leq C (1+ \vartheta_0(t,\,\xi)),
		\end{equation*}
		which gives
		\begin{equation*}
			|\sigma(\I B_2) + \sigma((\I B_2)^\ast)| \leq C(1+ \vartheta_0(t,\,\xi)).
		\end{equation*}
		Thus, we may conclude that we can estimate
		\begin{equation*}
			|\sigma(\I\psi(D_x)(A_1 + B_2+ R_2 + R_\infty) + \I\psi(D_x)(A_1 + B_2+ R_6)^\ast)| \leq C(1+ \vartheta_0(t,\,\xi)).
		\end{equation*}
		Choosing $K$ sufficiently large then gives the second statement.
	\end{Proof}

	\begin{Proposition}[Lemma 4.7 in \cite{Kinoshita.2005}]\label{Prop:Q1}
		For $s > 1$, we have that
		\begin{enumerate}[label = (\roman*)]
			\item $\psi(D_x) Q_1 \in L^\infty([0,\,T];\,\Zyg^s\OPS^0_{1,\,0})$,
			\item $\sigma\Big(\psi(D_x) \frac{Q_1 + Q_1^\ast}{2}\Big) \geq C I$.
		\end{enumerate}
	\end{Proposition}
	\begin{Proof}
		Let us write
		\begin{equation*}
			Z:= \bigg[\exp\bigg(- \int\limits_0^t \vartheta(s,\,D_x) \rmd s\bigg), A_1+B_2+R_2 + R_\infty\bigg] \exp\bigg( \int\limits_0^t \vartheta(s,\,D_x) \rmd s\bigg).
		\end{equation*}
		Applying Proposition~\ref{Prop:ZygComp} and Corollary~\ref{Cor:ZygComp} yields that
		\begin{align*}
			\sigma(Z) &= \sum\limits_{|\alpha| \leq 1} \partial_\xi^\alpha \sigma\bigg(\bigg[\exp\bigg(- \int\limits_0^t \vartheta(s,\,D_x) \rmd s\bigg), A_1+B_2+R_2 + R_\infty\bigg]\bigg) D_x^\alpha \exp\bigg( \int\limits_0^t \vartheta(s,\,\xi) \rmd s\bigg) + r\\
			&=\sum\limits_{|\alpha| = 1} \partial_\xi^\alpha \exp\bigg(- \int\limits_0^t \vartheta(s,\,\xi) \rmd s\bigg) D_x^\alpha (A_1+B_2+R_2 + R_\infty) \exp\bigg( \int\limits_0^t \vartheta(s,\,\xi) \rmd s\bigg) +r.
		\end{align*}
		Analyzing the symbols and applying Proposition~\ref{Prop:Theta} gives that $Z \in L^\infty([0,\,T];\,\Zyg^s\OPS^0_{1,\,0})$.
		
		Similarly to the argument in the proof of Proposition~\ref{Prop:Q0}, the second result follows for sufficiently large $K$ from the rules for computing the adjoint and the fact that $Z \in L^\infty([0,\,T];\,\Zyg^s\OPS^0_{1,\,0})$.
	\end{Proof}
	
	\subsection{Conclusion}\label{sec:conclusion}
	
	We consider the Cauchy problem
	\begin{equation*}
		\partial_t V + Q_0 V + Q_1 V = 0,\quad V(0,\,x) = V_0(x).
	\end{equation*}
	We have
	\begin{equation*}
		\partial_t\|V(t,\,\cdot)\|_{L^2}^2 = -2\Re(Q_0V,\,V) - 2\Re(Q_1V,\,V).
	\end{equation*}
	According to Propositions~\ref{Prop:Q0} and \ref{Prop:Q1} we have
	\begin{equation*}
		\sigma\Big(\psi(D_x) \frac{Q_0 + Q_0^\ast}{2}\Big) \geq 0, \qquad \sigma\Big(\psi(D_x) \frac{Q_1 + Q_1^\ast}{2}\Big) \geq 0,
	\end{equation*}
	which allows us to apply sharp G{\aa}rding's inequality (see Proposition~\ref{Prop:ZygGaar}).
	Since
	\begin{equation*}
		\psi(D_x) Q_1 \in L^\infty([0,\,T];\,\Zyg^s\OPS^0_{1,\,0}),
	\end{equation*}
	and 
	\begin{equation*}
		\psi(D_x) Q_0 \in L^\infty([0,\,T];\,\Zyg^s\OPS^{m_0}_{1,\,0}),
	\end{equation*}
	application of our result for sharp G{\aa}rding's inequal is only possible if \begin{equation*}
		s \geq \frac{2m_0}{2-m_0}.
	\end{equation*}
	This yields
	\begin{equation*}
		-2\Re(Q_0V,\,V) - 2\Re(Q_1V,\,V) \leq C \|V(t,\,\cdot)\|_{L^2}^2,
	\end{equation*}
	for $s \geq \frac{2m_0}{2-m_0}$.
	Application of Gronwall's lemma gives
	\begin{equation*}
		\|V(t,\,\cdot)\|_{L^2} \leq C_T \|V(t,\,\cdot)\|_{L^2}.
	\end{equation*}
	Using the definition of $U$ we obtain
	\begin{equation*}
		\| (\jbl D_x \jbr^{m-1} u(t,\,\cdot),\,\ldots,\,D_t^{m-1}u(t,\,\cdot))^T\|_{L^2}\\\leq C_T \| (\jbl D_x \jbr^{m-1} g_1(\cdot),\,\ldots,\,D_t^{m-1}g_m(\cdot))^T\|_{L^2},
	\end{equation*}
	where we used that $M_1,\,M_2,\,M_3$ and $ \exp\bigg(- \int\limits_0^t \vartheta(s,\,D_x) \rmd s\bigg)$ are operators of order zero.
	
	We conclude the proof by noting that condition~\eqref{cond:s} is due to the restriction $s \geq \frac{2m_0}{2-m_0}$ coming from sharp G{\aa}rding's inequality and due to the restriction $s > 1$ coming from our results for composition and adjoints in $\Zyg^s\OPS^m$.

	\section{Concluding remarks}\label{CONCLUSION}
	In this paper we consider the strictly hyperbolic Cauchy problem
	\begin{equation*}
	\begin{cases}
	D_t^m u - \sum\limits_{j = 0}^{m-1} \sum\limits_{|\gamma|+j = m} a_{m-j,\,\gamma}(t,\,x) D_x^\gamma D_t^j u = 0,\\
	D_t^{k-1}u(0,\,x) = g_k(x),\,k = 1,\,\ldots,\,m,
	\end{cases}
	\end{equation*}
	for $(t,\,x) \in [0,\,T]\times \R^n$, with coefficients low regular in time and space. By this we mean that the coefficients are less regular than Lipschitz in time and belong to the Zygmund space $\Zyg^s$ in $x$. Under suitable assumptions we prove a global (on $[0,\,T]$) well-posedness result without loss of derivatives if
	\begin{equation}\label{Conclusion:eq1}
	s \geq \max\{1+\ve,\,\frac{2m_0}{2-m_0}\},
	\end{equation}
	where $\ve > 0$ and $m_0\in(0,\,1]$ is related to the regular of the coefficients in time. The number $m_0$ is closer to $1$ for less regular (in time) coefficients and closer to $0$ if the coefficients are close to Lipschitz in time.
	
	Of course one would expect that a higher regular of the coefficients in time would allow us to lower the regularity of the coefficients with respect to $x$. Condition \eqref{Conclusion:eq1} however tells us that the index of the Zygmund spaces $s$ always has to be strictly larger than $1$. At some point, additional regularity in time of to coefficients does not improve the situation for $s$ any further. 
	
	\section*{Acknowledgements}
	The author wants to express his gratitude to Michael Reissig for many fruitful discussions and suggestions. Furthermore, he wants to thank Daniele Del Santo for his hospitality and the suggested improvements during the authors stay at Trieste University.

	\end{document}